\documentclass[review]{elsarticle}
\usepackage{setspace}
\usepackage{lineno,hyperref}
\usepackage{amsmath, amssymb,amsthm, mathtools, empheq}
\usepackage{paralist, mathtools}
\usepackage{caption}
\usepackage{subcaption}
\usepackage{graphics}
\usepackage{graphicx}
\usepackage{xcolor}
\usepackage{ulem}
\usepackage[title]{appendix}
\usepackage{placeins}
\usepackage{enumerate, comment}
\usepackage{mathrsfs}
\usepackage{mathtools}

\allowdisplaybreaks

\newtheorem{theorem}{Theorem}[section]

\newtheorem{lemma}[theorem]{Lemma}

\theoremstyle{definition}
\newtheorem{definition}[theorem]{Definition}
\newtheorem{remark}{Remark}

\newcommand{\ve}{\varepsilon}

\newcommand{\di}[1]{{\mathrm{div}}\left( #1\right)}
\newcommand{\eff}{{\mathrm{eff}}}
\newcommand{\tu}{\widetilde u_\ve}

\journal{}

\begin{document}

\begin{frontmatter}

\title{Derivation of a bidomain model for \\ bundles of myelinated axons}
\author[ad1]{Carlos Jerez-Hanckes}
\author[ad2]{Isabel A. Mart\'{i}nez \'{A}vila}
\author[ad3]{Irina Pettersson\corref{mycorrespondingauthor}}
\cortext[mycorrespondingauthor]{Corresponding author}
\ead{irinap@chalmers.se}
\author[ad4,ad3]{Volodymyr Rybalko}

\address[ad1]{Universidad Adolfo Ib\'a\~nez, Chile}
\address[ad2]{Pontificia Universidad Cat{\'o}lica de Chile, Chile}
\address[ad3]{Chalmers University of Technology and Gothenburg University, Sweden}
\address[ad4]{Institute for Low Temperature Physics and Engineering, Ukraine}

\begin{abstract}
The work concerns the multiscale modeling of a nerve fascicle of myelinated axons. We present a rigorous derivation of a macroscopic bidomain model describing the behavior of the electric potential in the fascicle based on the FitzHugh-Nagumo membrane dynamics. The approach is based on the two-scale convergence machinery combined with the method of monotone operators.      
\end{abstract}

\begin{keyword}
Nerve fascicle, myelinated axons, bidomain model, FitzHugh-Nagumo model, multiscale analysis, degenerate evolution equation.

\vspace{2mm}
\MSC[2010] Primary: 35K61, 35K57, 35B27; Secondary: 35K65, 92C30.

\end{keyword}

\end{frontmatter}

\linenumbers


\nolinenumbers

\section{Introduction}
 Modeling {the} electrical stimulation of nerves requires {biophysically consistent descriptions amenable also for computational purposes.} A typical nerve in the  
 peripheral nervous system contains several grouped fascicles, each of them comprising hundreds of axons \cite{standring2021gray}.
 This complex microstructure of neural  tissue presents an obvious problem for those attempting to desribe its 
 macroscopic response to electrical excitation. Specifically, one needs to know 
 both how signals propagate along a single axon and how axons influence each other in a bundle. 
 
 Electric currents along individual axons are usually modeled via cable theory, which dates back to works of W. Thomson (Lord Kelvin).
Fundamental insights into nerve cell excitability were made by A. Hodgkin and A. Huxley, who proposed a model that describes ionic mechanisms   underlying the initiation and propagation of action potentials in axons \cite{hodgkin1952}. 
Later a more simple model for nonlinear dynamics in axons was introduced in \cite{fitzhugh1955mathematical}, {known as the FitzHugh–Nagumo model}.

Multiscale homogenization techniques {were} used in recent works \cite{jerez2020derivation,jerez2021multiscale}  
to derive an effective cable equation describing propagation of signals in myelinated axons.
Ideas of homogenization theory can also be naturally applied {to account for} ephaptic coupling in bundles of axons, where 
neighboring axons can communicate via current flow through  the extracellular space. In 1978, experiments on giant squid axons 
were conducted \cite{ramon1978ephaptic} 
revealing evidence of ephaptic events and their physiological importance. 
%
%
Ephaptic interactions might be modelled by coupled systems of a large number of cable equations (as, e.g., in \cite{bokil2001ephaptic}, \cite{binczak2001ephaptic}), but a continuous mathematical model for a fascicle of myelinated axons, to our best knowledge, has not been rigorously derived. An analogous phenomenon of coupling 
is observed in the electrical conductance of cardiac tissues \cite{lin2010modeling}, leading to the celebrated {\it bidomain model}.  {It was first derived}
by J. Neu and W. Krassowska \cite{NEK93}. In \cite{franzone2002degenerate} the authors study the well-posedness of the reaction-diffusion systems modeling cardiac electric activity at the micro- and macroscopic level. They focus on the FitzHugh-Nagumo model (with recovery variable), and present a formal derivation of the effective bidomain model. The homogenization procedure is justified in \cite{pennacchio2005} where $\Gamma$-convergence is used for asymptotic analysis. Homogenization techniques based on two-scale convergence and unfolding
are applied in, e.g., \cite{collin2018mathematical}, \cite{BenMroSaaTal2019}, \cite{GraKar2019}, \cite{AmaAndTim2021} for  modeling of syncytial tissues.

The multiscale analysis of syncytial tissues includes the well-posedness of the microscopic problem, the homogenization procedure, and the well-posedness of the effective bidomain model. The latter question is interesting by itself, with solvability  {proven} using different approaches depending on the nonlinearity.
The solvability for a bidomain model in \cite{franzone2002degenerate} is based on a reformulation as a Cauchy problem for a variational evolution inequality in a properly chosen Sobolev space. This approach applies to the case of the FitzHugh-Nagumo equations. In \cite{Veneroni} existence and uniqueness {are given for solutions of} a wide class of models, including the classical Hodgkin-Huxley model, the first membrane model for ionic currents in an axon, and the Phase-I Luo-Rudy (LR1) model. In \cite{bourgault2009existence} the authors reformulate the coupled parabolic and elliptic PDEs into a single parabolic PDE by the introduction of a bidomain operator, which is a non-differential and non-local operator. This approach applies to fairly general ionic models, {such} as the Aliev-Panfilov and MacCulloch models.

  

 The asymptotic analysis of a nerve fascicle with a large number of axons also leads to a bidomain model.
 In \cite{mandonnet2011role} a linear model {is considered} without recovery variables. {Therein, it is} hypothesized that the homogenization procedure in \cite{pennacchio2005} leading to a macroscopic bidomain model for syncytical tissues can also be carried out for a fascicle of unmyelinated axons. 
 We extend this result to a nonlinear case and rigorously derive a bidomain model for a fascicle of myelinated axons. {In particular,} we consider {the} propagation of signals in a fascicle formed by a large number of 
 axons. The microstructure of the fascicle is {depicted as} a set of closely packed thin
cylinders---axons---with myelin sheaths arranged periodically in the surrounding  extracellular matrix. The characteristic microscale of the structure is given by a small parameter $\ve>0$. 
Distances between neighboring axons, their diameters and {the} spacing of unmyelinated parts of the axon's membrane---Ranvier nodes---are assumed {to be of} order $\ve$. By means of two-scale analysis we derive a bidomain model that
describes {the} asymptotic behavior of the transmembrane potential 
on Ranvier nodes when $\ve$ is sufficiently small.  We adopt the FitzHugh-Nagumo dynamics on the unmyelinated membrane. {Main} technical difficulties come from the nonlinear dynamics and {the}  lack of a priori estimates ensuring strong convergence of the membrane potential on the Ranvier nodes. This lack of compactness is caused by 
the fact that the axons form a disconnected microstructure inside the fascicle, which stands in the contrast with connected microstructure of syncytial tissues. 
In order to derive the homogenized problem we transform problem to a form allowing us to combine  two-scale convergence machinery with the method of monotone operators.
Well-posedness of the micro- and macroscopic problems are also shown via reduction to parabolic equations with monotone operators.

\section{Microscopic model}
\subsection{Problem setup}
\begin{figure}[htp]
\centering 
\def\svgwidth{.9\textwidth}
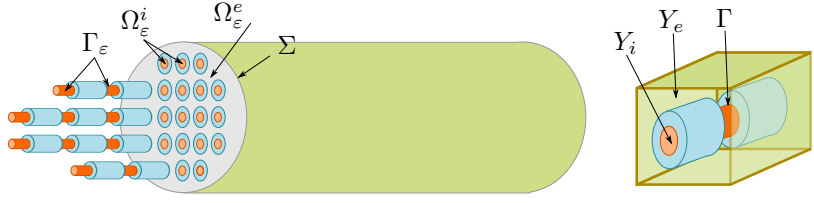
\caption{A fascicle of myelinated axons and the periodicity cell $Y$.}
\label{fig:bundle}
\end{figure}

A nerve fascicle is modeled by the cylinder   $\Omega:=(0,L)\times \omega\subset\mathbb{R}^3$ with length $L>0$ and cross section $\omega\subset \mathbb{R}^2$, being a bounded domain  in $\mathbb{R}^2$ with a Lipschitz boundary $\partial \omega$ {(see Figure \ref{fig:bundle})}. The lateral boundary of the cylinder is denoted by $\Sigma:=[0,L]\times \partial \omega$, with bases $S_0:=\{0\} \times \omega$, $S_L:=\{L\} \times \omega$.
The bulk of the cylinder consists of  an intracellular part formed by thin cylinders (axons), an extracellular part, and myelin sheaths. To describe the microstructure of the fascicle, we introduce a periodicity cell $Y:=[-\frac{1}{2},\frac{1}{2})\times [-R_0,R_0)^2$, consisting of three disjoint Lipschitz domains: (i) an intracellular part $Y_i:=[-\frac{1}{2},\frac{1}{2})\times D_{r_0}$, where $D_{r_0}$ is the disk with radius $0<r_0<\frac{1}{2}$; (ii) a myelin sheath $Y_m$; (iii) an extracellular domain $Y_e$. {The real positive radii satisfy $r_0<R_0$.} We denote by $\Gamma_{mi}:=\overline{Y_i} \cap \overline{Y_m}$ the interface between $Y_i$ and $Y_m$. The interface between the extracellular domain $Y_e$ and a myelin sheath $Y_m$ is $\Gamma_{me}:=\overline{Y_e} \cap \overline{Y_m}$.  The unmyelinated part of the boundary of $Y_i$ (the Ranvier node) will be denoted by $\Gamma=\overline{Y_i} \cap \overline{Y_e}$ (see Figure \ref{fig:bundle}). We will assume that $\Gamma$ does not degenerate, and, for simplicity, that $\Gamma$ is connected. 

The periodicity cell is translated by vertices of the lattice $\mathbb Z\times (2R_0 \mathbb Z)^2$ to form a $Y$-periodic structure, and then scaled by a small parameter $\ve >0$.  We take only those axons that are entirely contained in $\Omega$. As a result, the domain is the union of three disjoint parts $\Omega_\ve^i,  \Omega_\ve^e, \Omega_\ve^m$, and their boundaries (see Figure \ref{fig:bundle}). The unmyelinated part of the boundary of $\Omega_\ve^i$ is denoted by $\Gamma_\ve$. The 
boundary of the myelin is denoted by $\Gamma_\ve^m$.

Let $u_\ve$ denotes the electric potential $u_\ve = u_\ve^l$ in $\Omega_\ve^l$, $l=i, e$. We assume that $u_\ve$ satisfies homogeneous Neumann boundary conditions on the boundary of the myelin sheath $\Gamma_\ve^m$, i.e the myelin sheath is assumed to be a perfect insulator {(see \cite{jerez2020derivation} for other insulation assumptions)}.
The transmembrane potential $v_\ve = [u_\ve] = u_\ve^i - u_\ve^e$ is the potential jump across the Ranvier nodes $\Gamma_\ve$. We assume that the conductivity is a piecewise constant function:
\begin{align*}
a_\ve = \left\{
\begin{array}{l}
a_e \quad \mbox{in}\,\, \Omega_\ve^e,\\
a_i \quad \mbox{in}\,\, \Omega_\ve^i.
\end{array}
\right.
\end{align*}
On $\Gamma_\ve$ we {further} assume current continuity, and FitzHugh-Nagumo \cite{fitzhugh1955mathematical,nagumo1962active} dynamics for the transmembrane potential.
Namely,  the ionic current is {described as}
\begin{align*}
I_{ion}(v_{\ve}, g_{\ve}) = \frac{v^3_{\ve}}{3} - v_{\ve} - g_{\ve},
\end{align*}
where $g_{\ve}$ is the recovery variable whose evolution is governed by the ordinary differential equation
\begin{align*}
    \partial_t {g}_{\ve} = \theta v_{\ve} + a - bg_{\ve}
\end{align*}
with constant coefficients $\theta, a, b >0$. 
{Thus, the} electric activity in the bundle {$\Omega$} is described by the following system of equations for the unknowns  $v_\ve$ and $g_\ve$:
\begin{align}
&-\di{a_\ve \nabla u_\ve}  =0, \, &(t,x) &\in (0,T)\times (\Omega_\ve^i \cup \Omega_\ve^e), \nonumber \\
&a_e \nabla u_\ve^e\cdot \nu  = a_i \nabla u_\ve^i\cdot \nu, \, &(t,x) &\in (0,T)\times \Gamma_\ve, \nonumber\\
\nonumber
&\ve (c_m \partial_t[u_\ve] + I_{ion}([u_\ve], g_\ve))
= - a_i \nabla u_\ve^i\cdot \nu
, \, &(t,x) & \in(0,T)\times  \Gamma_\ve,
\\
&\partial_t g_\ve  = \theta [u_\ve] + a - bg_\ve, \, &(t,x) & \in(0,T)\times  \Gamma_\ve, \nonumber \\
&u_\ve   = 0, \, &(t,x)& \in (0,T)\times (S_{0} \cup S_L), \label{eq:orig-prob}
\\
&a_e \nabla u_\ve^e\cdot \nu   = J_\ve^e(t,x), \, &(t,x)& \in (0,T)\times\Sigma, \nonumber \\
&\nabla u_\ve^e\cdot \nu   = 0, \, &(t,x)& \in (0,T)\times\Gamma_\ve^m, \nonumber 
\\
& [u_\ve](0,x) =V_\ve^0(x), \,\,g_\ve(0,x)=G_\ve^0(x), \, &x& \in \Gamma_\ve, \nonumber
\end{align}
where $\nu$ denotes the unit normal on $\Gamma_\ve$, $\Gamma_\ve^m$, and $\Sigma$, exterior to $\Omega_\ve^i$, $\Omega_\ve^m$, and $\Omega$, respectively. 
The function $J_\ve^e(t,x)$ models an external {boundary} excitation of the nerve fascicle. 

We study the asymptotic behavior of $u_\ve$, as $\ve \to 0$, and derive a macroscopic model describing the potential $u_\ve$ in the fascicle, under the following conditions:
\begin{itemize}
\item[(H1)]
The initial data is such that\footnote{Throughout, $C$ denotes a generic constant independent of $\ve$, whose value may be different from line to line.} $\|V_\ve^0\|_{L^4(\Gamma_\ve)} \le C$. Moreover, we assume that $V_\ve^0= V_\ve^i - V_\ve^e$, where $V_\ve^l$, $l=i, e$, can be extended to the whole $\Omega$ such that, keeping the same notation for the extension, $\|V_\ve^l\|_{H^1(\Omega)} \le C$ and $V_\ve^l=0$ on $S_0\cup S_L$. We also assume that there exist weak limits of $V_\ve^l \rightharpoonup  V^l$ in $H^1(\Omega)$, $l=i, e$.\\[2mm]

\item[(H2)]
There exists $G^0\in L^2(\Omega)$,  such that
\begin{itemize}
    \item[$\bullet$] for any $\phi \in C(\overline \Omega)$, {it holds that}
\begin{align*}
\lim\limits_{\ve \to 0}\ve  \int_{\Gamma_\ve} G_\ve^0(x) \phi(x)\, d\sigma =
\frac{|\Gamma|}{|Y|}\int_\Omega G^0(x) \phi(x) \,dx;
\end{align*}
\item[$\bullet$]
$
\displaystyle
\ve \int_{\Gamma_\ve} |G_\ve^0|^2\, d\sigma \to \frac{|\Gamma|}{|Y|} \int_\Omega |G^0|^2\, dx, \quad \ve \to 0.
$
\end{itemize}

\vspace{2mm}
\item[(H3)]
The external excitation $J_\ve^e \in L^2((0,T)\times \Sigma)$ converges weakly to $J^e(t,x)$, as $\ve \to 0$, and
\begin{align*}
\int_0^t \int_\Sigma |\partial_\tau J_\ve^e |^2\, d\sigma d\tau \le C.
\end{align*}
\end{itemize}


\begin{remark} Hypothesis (H2) actually assumes strong two-scale convergence (cf. Proposition 2.5 in \cite{AlDa-95}).
Note that the hypothesis (H2) is satisfied if $G_\ve^0$ is sufficiently regular, e.g. continuous, and independently of $\ve$. 
\end{remark}

\subsection{Main result}
The main result of the paper is given in short form by Theorem \ref{th:main-short}, showing that the asymptotic behavior of solutions of {the boundary value problem} \eqref{eq:orig-prob} is described by the following effective bidomain model in $\Omega$:
\begin{align}
& c_m \partial_t v_0 + I_{ion}(v_0, g_0)  = a_i^\eff \partial_{x_1 x_1}^2 u_0^i,  &(t,x) &\in (0,T)\times \Omega,\nonumber\\
& c_m \partial_t v_0 + I_{ion}(v_0, g_0)  = -\di{a_e^\eff \nabla u_0^e},  &(t,x) &\in (0,T)\times \Omega,\nonumber\\
\label{eq:hom-prob}
&\partial_t g_0  = \theta v_0+a - b\, g_0,  &(t,x) &\in (0,T)\times \Omega,\\
&u_0^{i,e}(t,x) =0,&(t,x) &\in (0,T)\times(S_0\cup S_L),\nonumber\\
&a_e^\eff \nabla u_0^{e} \cdot \nu = J^e, &(t,x) &\in (0,T)\times \Sigma,\nonumber\\
&v_0(0,x) =V^i(x) - V^e(x), \,\, g_0(0,x)=G^0(x),
\hskip -0.5cm & x &\in \Omega, \nonumber
\end{align}
where $v_0=u_0^i - u_0^e$.
The effective scalar coefficient $a_i^\eff$ is
\begin{align}
\label{eq:eff-coef-i}
a_i^\eff:= \frac{|Y_i|}{|\Gamma|}a_i.
\end{align}
The effective matrix $a_{e}^\eff {\in\mathbb{R}^{3\times3}}$ is given by
\begin{align}
\label{eq:eff-coef-e}
(a_e^\eff)_{kl} := \frac{1}{|\Gamma|} \int_{Y_e} a_e (\partial_{l} N_k^e(y) +\delta_{kl})\, dy,\quad {k,l=1,2,3,}
\end{align}
with the functions $N_k^e$, $k=1, 2, 3$, solving the following auxiliary cell problems in $Y_e$
\begin{align*}
-&\Delta N_k^e = 0, &y&\in Y_e, \nonumber\\
&\nabla N_k^e \cdot \nu = -\nu_k, \hskip -2.5cm&y&\in \Gamma\cup \Gamma_m,
\\
&N_k^e(y) \,\,\,\mbox{is} \,\,\, Y-\mbox{periodic}.\hskip -2.5cm&& \nonumber
\end{align*}
\begin{theorem}
\label{th:main-short}
Under the hypothesis (H1)--(H3), the solutions $v_\ve=[u_\ve]$, $g_\ve$ of the microscopic problem \eqref{eq:orig-prob} converge to the solution $v_0=u_0^i-u_0^e$, $g_0$ of the macroscopic {one} \eqref{eq:hom-prob} in the following sense:
\begin{itemize}
\item[(i)] For any $\Phi(t,x)\in C([0,T]\times \overline{\Omega})$, {it holds that}
\begin{align*}
\lim_{\ve\to 0} \ve \int_0^T \int_{\Gamma_\ve} v_\ve(t, x) \Phi(t, x)\, d\sigma_{x}dt
= \frac{|\Gamma|}{|Y|} \int_0^T \int_{\Omega} v_0(t, x) \Phi(t, x)\, dx dt.
\end{align*}
\item[(ii)] For any $t\in [0,T]$, {one has} $$ \lim_{\ve \to 0}\ve \int_{\Gamma_\ve}|v_\ve|^2\, d\sigma = \frac{|\Gamma|}{|Y|} \int_{\Omega} |v_0|^2\, dx.$$

\item[(iii)]
For any $\Phi(t,x)\in C([0,T]\times \overline{\Omega})$,
\begin{align*}
\lim_{\ve\to 0} \ve \int_0^T \int_{\Gamma_\ve} g_\ve(t, x) \Phi(t, x)\, d\sigma_{x} dt
= \frac{|\Gamma|}{|Y|} \int_0^T \int_{\Omega} g_0(t, x) \Phi(t, x)\, dx dt.
\end{align*}
\item[(iv)]
For any $t\in [0,T]$, $\displaystyle \ve \int_{\Gamma_\ve}|g_\ve|^2\, d\sigma \to \frac{|\Gamma|}{|Y|} \int_{\Omega} |g_0|^2\, dx$, as $\ve \to 0$.

\item[(v)]
$\displaystyle \int_0^T \int_{\Omega_\ve^{i,e}} |u_\ve^{i,e} - u_0^{i,e}|^2\, dx dt \to 0$, as $\ve \to 0$.
\end{itemize}
\end{theorem}
\begin{remark}
If $v_0$ is continuous, the convergence (i), (ii) implies strong convergence of $v_\ve$. Namely, for any $t\in [0,T]$, {one obtains}
\begin{align*}
\lim_{\ve \to 0}\, \ve \int_{\Gamma_\ve}|v_\ve-v_0|^2\, d\sigma= 0.
\end{align*}
In general, approximating $v_0$ in $L^2(\Omega)$ by $v_{0\delta}\in C(\Omega)$, we have
\begin{align*}
\limsup_{\delta \to 0} \limsup_{\ve \to 0}\ve \int_{\Gamma_\ve}|v_\ve-v_{0\delta}|^2\, d\sigma= 0.
\end{align*}
\end{remark}
\begin{remark}
The result can be generalized to the case of {a} varying cross section, as in \cite{jerez2021multiscale}. In such case, the solution $N_1^i$ of the cell problem \eqref{eq:cell-prob-ii} is no longer constant, and the corresponding effective coefficient is given by
\begin{align*}
    a_i^\eff = \frac{1}{|\Gamma|}\int_{Y_i} a_i (\partial_1 N_1^i + 1)dy.
\end{align*}
\end{remark}

\begin{remark}
Hypothesis (H2) can be generalized to the case of {an} oscillating initial function $G_\ve^0$. Namely, assume that there exists $G^0(x,y) \in L^2(\Omega \times \Gamma)$, $Y$-periodic in $y$ such that
\begin{itemize}
    \item[$\bullet$] for any $\phi(x,y)\in C(\overline \Omega\times Y)$, $Y$-periodic in $y$,
\begin{align*}
\lim\limits_{\ve \to 0}\ve  \int_{\Gamma_\ve} G_\ve^0(x) \phi\left(x,\frac{x}{\ve}\right)\, d\sigma_{x} =
\frac{1}{|Y|}\int_\Omega \int_\Gamma G^0(x,y) \phi(x,y) \,d\sigma_y dx;
\end{align*}
\item[$\bullet$]
$
\displaystyle
\ve \int_{\Gamma_\ve} |G_\ve^0|^2\, d\sigma \to \frac{1}{|Y|} \int_\Omega \int_\Gamma |G^0(x,y)|^2\, d\sigma_y dx, \quad \ve \to 0.$
\end{itemize}
Then, the two-scale limit $\widetilde{g}_0(t,x,y)$ of $g_\ve$ does depend on the fast variable $y$, and denoting $\displaystyle g_0(t,x)= \frac{1}{|\Gamma|} \int_\Gamma \widetilde{g}_0(t,x,y)\, d\sigma_y$, the effective problem reads
\begin{align*}
& c_m \partial_t v_0 + I_{ion}(v_0, g_0)  = a_i^\eff \partial_{x_1 x_1}^2 u_0^i,  &(t,x) &\in (0,T)\times \Omega,\nonumber\\
& c_m \partial_t v_0 + I_{ion}(v_0, g_0)  = -\di{a_e^\eff \nabla u_0^e},  &(t,x) &\in (0,T)\times \Omega,\nonumber\\
&\partial_t \widetilde{g}_0  = \theta v_0+a - b\, \widetilde{g}_0,  &(t,x,y) &\in (0,T)\times \Omega\times Y,\\
&u_0^{i,e}(t,x) =0,&(t,x) &\in (0,T)\times(S_0\cup S_L),\nonumber\\
&a_e^\eff \nabla u_0^{e} \cdot \nu = J^e, &(t,x) &\in (0,T)\times \Sigma,\nonumber\\
&v_0(0,x) =V^i(x) - V^e(x), \,\, \widetilde{g}_0(0,x)=G^0(x,y)
\hskip -0.5cm & x &\in \Omega, \,\, y\in Y. \nonumber
\end{align*}
Thanks to the linearity of the equation $\partial_t \widetilde{g}_0  = \theta v_0+a - b\, \widetilde{g}_0$, averaging in $y$, yields \eqref{eq:hom-prob} with the initial condition  $\displaystyle g_0(0,x)=\frac{1}{|\Gamma|}\int_{\Gamma} G^0(x,y)\, d\sigma_y$.
\end{remark}

\vspace{5mm}
\subsection{Well-posedness}
\label{sec:existence}
{In order to} show the well-posedness of the microscopic problem \eqref{eq:orig-prob}, we write it as a Cauchy problem for an abstract parabolic equation.

We multiply \eqref{eq:orig-prob} by {a} smooth function $\phi=\begin{cases}\phi^{i}\,\, \mbox{in}\,\, \Omega_\ve^{i}\\\phi^{e}\,\, \mbox{in}\,\, \Omega_\ve^{e}
\end{cases}$, $\phi^{i,e}=0$ on $S_0\cup S_L$, and integrate by parts:
\begin{align*}
&\ve \int_{\Gamma_\ve} c_m \partial_t v_\ve [\phi] \, d\sigma
+ \int_{\Omega_\ve^i \cup \Omega_\ve^e}a_\ve \nabla u_\ve \cdot \nabla \phi \, dx
+ \ve \int_{\Gamma_\ve} I_{ion}(v_\ve, g_\ve) [\phi] \, d\sigma
=
 \int_{\Sigma} J_\ve^e \phi \, d\sigma.
\end{align*}
Let us introduce an auxiliary function $q_\ve$ solving the following problem:
\begin{align}
\label{eq:stat-prob-2}
-&\di{a_\ve \nabla  q_\ve}  =0, \, &x & \in \Omega_\ve^i \cup \Omega_\ve^e \cup \Gamma_\ve, \nonumber\\
&\nabla q_\ve\cdot \nu   = 0, \, &x& \in \Gamma_{m,\ve},  \\
&a_e \nabla q_\ve\cdot \nu   = J_\ve^e(t,x), \, &x& \in \Sigma, \nonumber \\
& q_\ve   = 0, \, &x& \in (S_{0} \cup S_L). \nonumber
\end{align}
Since the jump of $q_\ve$ through the Ranvier nodes $\Gamma_\ve$ is zero, the change of unknown 
\begin{align*}
\tu=u_\ve - q_\ve
\end{align*}
allows us to transfer the external excitation $J_\ve^e$ from the lateral boundary $\Sigma$ to the membrane $\Gamma_\ve$. Namely, we get the following weak formulation for the new unknown function $\tu$: 
\begin{align*}
&\ve \int_{\Gamma_\ve} c_m \partial_t v_\ve [\phi] \, d\sigma
+ \int_{\Omega_\ve^i \cup \Omega_\ve^e}a_\ve \nabla \tu \cdot \nabla \phi \, dx
+ \ve \int_{\Gamma_\ve} I_{ion}(v_\ve, g_\ve) [\phi] \, d\sigma\\
&+ \int_{\Gamma_\ve} (a_i \nabla q_\ve\cdot \nu) [\phi]\, d\sigma
=0. \nonumber
\end{align*}
{Let us define the subspace 
$$ H^1_{S_0\cup S_L}(\Omega_\ve^i \cup \Omega_\ve^e):=\left\{\phi\in H^1_{S_0\cup S_L}(\Omega_\ve^i \cup \Omega_\ve^e) \ : \ \phi\big|_{S_0\cup S_L}=0\right\},$$
and} introduce {the} operator
$A_\ve: D(A_\ve) \subset H^{1/2}(\Gamma_\ve) \to H^{-1/2}(\Gamma_\ve)$ {as follows}
\begin{align}
\label{def:A_eps}
(A_\ve v_\ve, [\phi])_{L^2(\Gamma_\ve)}:= \int_{\Omega_\ve^i \cup \Omega_\ve^e}  a_\ve \nabla \tu \cdot \nabla \phi \, dx, \quad
\forall \ \phi\in H^1_{S_0\cup S_L}(\Omega_\ve^i \cup \Omega_\ve^e),
\end{align}
where $\tu \in H^1(\Omega_\ve^i \cup \Omega_\ve^e)$, for a given jump $[\tu]=v_\ve$, solves the following problem: 
\begin{align}
\label{eq:stat-prob-1}
&-\di{a_\ve \nabla  \tu}  =0, \, &x & \in \Omega_\ve^i \cup \Omega_\ve^e, \nonumber\\
&a_e \nabla \tu^e\cdot \nu  = a_i \nabla \tu^i\cdot \nu, \, &x &\in  \Gamma_\ve, \nonumber\\
&\tu^i - \tu^e=v_\ve, \, &x& \in  \Gamma_\ve, 
\\
&a_\ve \nabla \tu\cdot \nu   = 0, \, &x& \in \Gamma_{m,\ve}, \nonumber \\
&a_e \nabla \tu\cdot \nu   = 0, \, &x& \in \Sigma, \nonumber \\
& \tu   = 0, \, &x& \in (S_{0} \cup S_L). \nonumber
\end{align}
{Thus}, the problem can be rewritten in the following {compact} form:
\begin{align}
\label{eq:system-1}
  & \ve c_m \partial_t v_\ve+ A_\ve v_\ve
+ \ve I_{ion}(v_\ve, g_\ve)
=-a_i\nabla q_\ve \cdot \nu,\\
& \partial_t g_\ve + b g_\ve - \theta v_\ve =a
\nonumber
\end{align}
on $\Gamma_\ve$.
In order to reduce the problem to a monotone {one}, we {perform} the following change of unknowns: 
\begin{align}
\label{def:W_eps}
W_\ve =\begin{pmatrix} w_\ve\\[2mm] h_\ve \end{pmatrix}
=e^{-\lambda t}\begin{pmatrix} v_\ve\\[2mm] g_\ve\end{pmatrix}, \quad W_\ve^0 = \begin{pmatrix} V_\ve^0\\[2mm] G_\ve^0 \end{pmatrix}.
\end{align}
{with $\lambda$ real positive.} Substituting \eqref{def:W_eps} into \eqref{eq:system-1} yields
\begin{align*}
\ve \partial_t \begin{pmatrix} w_\ve \\[6mm] h_\ve \end{pmatrix}
&+ 
\begin{pmatrix}
\displaystyle \frac{1}{c_m}A_\ve w_\ve + \frac{\ve}{c_m}\left(\frac{e^{2\lambda t}}{3}w_\ve^3 - w_\ve - h_\ve\right) + \ve \lambda w_\ve \\[6mm]
\ve (b+\lambda) h_\ve -\ve \theta w_\ve 
\end{pmatrix}\\
&=
e^{-\lambda t} \begin{pmatrix}
\displaystyle
- \frac{a_i}{c_m}\nabla q_\ve \cdot \nu\\[6mm]
\ve a 
\end{pmatrix},
\end{align*}
which can be further rewritten as follows:
\begin{align}
\label{eq:abstract-parabol}
&\ve \partial_t W_\ve + \mathbb A_\ve(t, W_\ve) = F_\ve(t), \quad (t,x) \in (0,T)\times \Gamma_\ve,\\
&W_\ve(0,x)=W_\ve^0(x), \quad x\in \Gamma_\ve. \nonumber
\end{align}
\begin{align}
\label{eq:A}
    \mathbb{A}_\ve(t, W_\ve) &:= B_\ve^{(1)}(t, W_\ve) + B_\ve^{(2)}(t, W_\ve),\\
    \label{eq:B1}
    B_\ve^{(1)}(t, W_\ve)&:=\begin{pmatrix}
\displaystyle \frac{1}{c_m}A_\ve w_\ve + \ve \left(\lambda -\frac{1}{c_m}\right)w_\ve - \frac{\ve}{c_m}h_\ve\\[4mm]
\ve (b+\lambda) h_\ve -\ve \theta w_\ve 
\end{pmatrix},\\
\label{eq:B2}
B_\ve^{(2)}(t, W_\ve) &:=
\begin{pmatrix}
\displaystyle
\ve \frac{e^{2\lambda t}}{3c_m} w_\ve^3 \\[4mm] 0
\end{pmatrix}, \quad
F_\ve(t) := e^{-\lambda t} \begin{pmatrix}
\displaystyle
-\frac{a_i}{c_m}\nabla q_\ve \cdot \nu \\[2mm]
\ve a 
\end{pmatrix}.
\end{align}
Here the operator $A_\ve$ is defined in \eqref{def:A_eps}.

The existence of a unique solution to problem \eqref{eq:abstract-parabol} follows from Theorem 1.4 in \cite{lions1969quelques} and Remark 1.8 in Chapter 2 (see also Theorem 4.1 in \cite{showalter2013monotone}). For the reader's convenience, we formulate the corresponding result below.
\begin{lemma}
\label{lm:our-existence}
Let $V_i$, $i=1,\ldots, m$, be reflexive Banach spaces, and $H$ be a real Hilbert space such that $V_i\subset H\subset V_i'$. Let $A(t)=\sum_{i=1}^m A_i(t)$, and let $\{A_i(t); \,\,t\in [0,T]\}$, $i=1, \ldots, m$, be a family of nonlinear, monotone, and demi-continuous operators from $V_i$ to $V_i'$ that satisfy the following conditions:
\begin{itemize}
    \item[(i)]
    The function $t \mapsto A_i(t)u(t)\in V_i'$ is measurable for every measurable function $u:[0,T] \to V$.
    \item[(ii)] There exists a seminorm $[u]$ on $V_i$ such that, for some constants $\alpha_1>0$ and $\alpha_2>0$, we have {that}
    \begin{align*}
        [u] + \alpha_1 \|u\|_H \ge \alpha_2 \|u\|_{V_i},
    \end{align*}
    and for some $\overline{c}>0$ and $p_i>1$,
    \begin{align*}
        (A_i(t)u,u)\ge \overline{c} [u]^{p_i}, \quad u \in V_i, \,\, t\in [0,T].
    \end{align*}
    \item[(iii)] For some $\underline{C}$ and the same $p_i>1$ as in (ii), 
    $$
    \|A_i(t)u\|_{V_i'} \le \underline{C}(1+\|u\|_{V_i}^{p_i-1}), \quad  u\in V_i, \,\,t\in [0,T].
    $$ 
\end{itemize}
Then, for every $u_0 \in H$ and $f \in \sum_{i=1}^m L^{q_i}(0,T; V_i')$, $1/p_i + 1/q_i=1$, there is a unique absolutely continuous function $u\in \cap_{i=1}^m W^{1,q_i}([0,T]; V_i')$ that satisfies
\begin{align*}
    & u\in L^\infty([0,T]; H), \,\, u\in \cap_{i=1}^m L^{p_i}([0,T]; V_i),\\
    &\frac{du}{dt}(t) + A(t) u(t) = f(t), \quad \mbox{a.e.} \,\, t\in (0,T),\\
    & u(0)=u_0.
\end{align*}
\end{lemma}
In order to apply Lemma \ref{lm:our-existence},  we introduce the necessary functional spaces:
\begin{align*}
&H= L^2(\Gamma_\ve) \times  L^2(\Gamma_\ve),\\
&\widetilde H^{1/2}(\Gamma_\ve)=\big\{v=(u^i-u^e)\Big|_{\Gamma_\ve}:\,\, u^l\in H^1(\Omega_\ve^l), \,\, u^l=0 \,\, \mbox{on}\,\, S_0\cap S_L, l=i, e \Big\},\\
&V_1 = \widetilde H^{1/2}(\Gamma_\ve)\times L^2(\Gamma_\ve), \quad V_1'=H^{-1/2}(\Gamma_\ve)\times L^2(\Gamma_\ve),\\
&V_2 = L^4(\Gamma_\ve) \times L^2(\Gamma_\ve), \quad V_2'= L^{4/3}(\Gamma_\ve) \times L^2(\Gamma_\ve).
\end{align*}

As the operator $A_1(t, \cdot): V_1 \to V_1'$ we take $B_\ve^{(1)}(t, \cdot)$ given by \eqref{eq:B1}; as the operator $A_2(t, \cdot): V_2 \to V_2'$ we take $B_\ve^{(2)}(t, \cdot)$ given by \eqref{eq:B2}.
Let us check that the operator $\mathbb A_\ve(t, \cdot)=B_\ve^{(1)} + B_\ve^{(2)}$ satisfies the assumptions of Lemma \ref{lm:our-existence}. The right-hand side $F_\ve$ satisfies clearly the assumptions of Lemma \ref{lm:our-existence}. 

\begin{lemma}
\label{lm:properties-B1}
For every $t\in[0,T]$, the linear operator $B_\ve^{(1)}(t,\cdot): V_1 \to V_1'$ has the following properties:
\begin{itemize}
    \item[(i)] Monotonicity:
    \begin{align*}
        (B_\ve^{(1)}(t,W_1) - B_\ve^{(1)}(t,W_2), W_1 - W_2) \ge 0, \quad \forall \ W_1, W_2 \in V_1.
    \end{align*}
    \item[(ii)] Coercivity:
    \begin{align*}
        (B_\ve^{(1)}(t,W), W) \ge C_1 \|W\|_{V_1}^2, \quad \forall \ W\in V_1.
    \end{align*}
    \item[(iii)] Boundedness: 
    \begin{align*}
        \|B_\ve^{(1)}(t,W)\|_{V_1'} \le C_2 \|W\|_{V_1}, \quad \forall\  W\in V_1.
    \end{align*}
\end{itemize}
\end{lemma}
\begin{proof} (i) The monotonicity of the operator $B_\ve^{(1)}$ follows from its linearity and coercivity properties (as shown below).

\noindent
(ii) By \eqref{eq:B1}, for any $W_\ve \in \widetilde H^{1/2}(\Gamma_\ve) \times L^2(\Gamma_\ve)$, we have
\begin{align*}
(B_\ve^{(1)}(t,W_\ve), W_\ve)
&= \frac{1}{ c_m} \int_{\Omega_\ve^i \cup \Omega_\ve^e} a_\ve |\nabla \widetilde w_\ve|^2\, dx
+ \ve \left(\lambda - \frac{1}{c_m}\right)\int_{\Gamma_\ve} |w_\ve|^2\, d\sigma\\
&- \ve \left(\theta +\frac{1}{c_m} \right)\int_{\Gamma_\ve} h_\ve w_\ve\, d\sigma
+ \ve(b+\lambda) \int_{\Gamma_\ve} |h_\ve|^2\, d\sigma.
\end{align*}
Here $\widetilde w_\ve= e^{-\lambda t} u_\ve$ solves  \eqref{eq:stat-prob-1} with the jump on $\Gamma_\ve$ that equals to $e^{-\lambda t} v_\ve$.
Using the trace inequality and choosing $\lambda$ sufficiently large and independent of $\ve$, we obtain
\begin{align*}
(B_\ve^{(1)}(t,W_\ve), W_\ve)
\ge C_1^\ve \|w_\ve\|_{\widetilde H^{1/2}(\Gamma_\ve)}^2
+ C_2^\ve \|h_\ve\|_{L^2(\Gamma_\ve)}^2 = C^\ve \|W_\ve\|_{V_1}^2.
\end{align*}
Here $C_1^\ve, C_2^\ve$, and $C^\ve$ are positive constants.\\
\noindent
(iii) Let us estimate the norm of $B_\ve^{(1)}(t, W)$. For any $W_\ve \in V_1$ and a test function $\Phi=([\varphi], \psi)^T \in V_1$, by \eqref{eq:A} we have 
\begin{align*}
(B_\ve^{(1)}(t, W_\ve), \Phi)_{L^2(\Gamma_\ve)^2}
&= \frac{1}{c_m} \int_{\Omega_\ve^i \cup \Omega_\ve^e} a_\ve \nabla \widetilde w_\ve \cdot \nabla \varphi\, dx
+\ve \left(\lambda - \frac{1}{c_m}\right) \int_{\Gamma_\ve} w_\ve [\varphi]d\sigma\\
&-\frac{\ve}{c_m}\int_{\Gamma_\ve} h_\ve [\varphi]d\sigma
+ \ve (b+\lambda)\int_{\Gamma_\ve} h_\ve \psi d\sigma
-\ve \theta \int_{\Gamma_\ve} w_\ve \psi d\sigma.
\end{align*}
Here $\varphi$ solves a stationary problem \eqref{eq:stat-prob-1} with a given jump $[\varphi]$ on $\Gamma_\ve$. 
Clearly, $\|\nabla \widetilde w_\ve\|_{L^2(\Omega_\ve^i \cup \Omega_\ve^e)} \le C \|w_\ve\|_{\widetilde H^{1/2}(\Gamma_\ve)}$. The test function $\varphi$ is estimated in a standard way in terms of $\|[\varphi]\|_{\widetilde H^{1/2}(\Gamma_\ve)}$. Then, by the Cauchy-Schwartz inequality, {one retrieves}
\begin{align*}
(B_\ve^{(1)}(t, W_\ve), \Phi)_{L^2(\Gamma_\ve)^2}
&\le C_1\|w_\ve\|_{\widetilde H^{1/2}(\Gamma_\ve)} \|[\varphi]\|_{\widetilde H^{1/2}(\Gamma_\ve)}\\
&+
C_2 (\|w_\ve\|_{\widetilde H^{1/2}(\Gamma_\ve)} + \|h_\ve\|_{\widetilde H^{1/2}(\Gamma_\ve)})
\|[\Phi]\|_{V_1},
\end{align*}
which proves the estimate from above for $\|B_\ve^{(1)}(t, W)\|_{V_1'}$.
\end{proof}

\vspace{5mm}
\begin{lemma}
\label{lm:properties-B2}
For every $t\in[0,T]$, the operator $B_\ve^{(2)}(t,\cdot): V_2 \to V_2'$ has the following properties: 
\begin{itemize}
    \item[(i)] Monotonicity:
    \begin{align*}
        (B_\ve^{(2)}(t,W_1) - B_\ve^{(2)}(t,W_2), W_1 - W_2) \ge 0, \quad \forall \ W_1, W_2 \in V_2.
    \end{align*}
    \item[(ii)] Coercivity: $\|\cdot\|_{L^4(\Gamma_\ve)}$
    defines a seminorm on $V_2$ such that, for some constants $\alpha_1>0$ and $\alpha_2>0$, we have  
    \begin{align*}
        \|W\|_{L^4(\Gamma_\ve)} + \alpha_1\|W\|_H \ge \alpha_2 \|W\|_{V_2},
    \end{align*}
    and
    \begin{align*}
        (B_\ve^{(2)}(t,W), W) \ge C_1 \|W\|_{V_2}^4, \quad \forall \ W\in V_1.
    \end{align*}
    \item[(iii)] Boundedness:
    \begin{align*}
        \|B_\ve^{(2)}(t,W)\|_{V_2'} \le C_2 \|W\|_{L^4(\Gamma_\ve)}^4, \quad \forall \ W\in V_2.
    \end{align*}
\end{itemize}
\end{lemma}

\begin{proof}
(i) The monotonicity of $B_\ve^{(2)}$ follows from the monotonicity of the cubic function $f(u)=u^3$. \\

\noindent
(ii) By definition \eqref{eq:B2}, {it holds that}
\begin{align*}
(B_\ve^{(2)}(t, W_\ve), W_\ve) = \frac{\ve e^{2\lambda t}}{3c_m} \int_{\Gamma_\ve} |w_\ve|^4\, d\sigma,
\end{align*}
which proves (ii).

\noindent
(iii) The boundedness follows from \eqref{eq:B2}:
\begin{align*}
\|B_\ve^{(2)}(t, W_\ve)\|_{V_2'} = 
\ve \left[\int_{\Gamma_\ve} \left(\frac{e^{2\lambda t}}{3c_m} (w_\ve)^3\right)^\frac{4}{3}\, d\sigma \right]^\frac{3}{4}
= \frac{\ve e^{2\lambda t}}{3c_m} \|w_\ve\|_{L^4(\Gamma_\ve)}^{3}
\le C^\ve \|W_\ve\|_{V_2}^3,
\end{align*}
where $C^\ve$ is a positive constant. 
\end{proof}
Obviously, the function $t\mapsto \mathbb A_\ve(t, W)$ satisfies the measurability assumption of Lemma \ref{lm:our-existence}, and the demi-continuity property follows from the estimates in Lemmas \ref{lm:properties-B1} and \ref{lm:properties-B2}.


\subsection{A priori estimates}

The next lemma provides the estimates for $(z_\ve, h_\ve)=e^{-\lambda t}(u_\ve, g_\ve)$, where $[z_\ve]=w_\ve$, at time $t=0$. 
\begin{lemma}
\label{lm:t=0}
Under hypotheses (H1)--(H3), at time $t=0$ the following estimates hold
\begin{align}
\label{est:t=0}
\int_{\Omega_\ve^i\cup \Omega_\ve^e} a_\ve |\nabla z_\ve|^2\, dx \Big|_{t=0} +
\int_{\Sigma} |z_\ve|^2\, d\sigma \Big|_{t=0} \le C.
\end{align}
\end{lemma}
\begin{proof}
One can see that the operator $A_\ve$ given by \eqref{eq:A} can be defined by means of the minimization problem
\begin{align*}
(A_\ve w_\ve, w_\ve) = \min_{[\phi_\ve]=w_\ve} 
\int_{\Omega_\ve^i\cup \Omega_\ve^e} a_\ve |\nabla \phi_\ve|^2\, dx,
\end{align*}
where the minimum is taken over the functions $\phi_\ve \in H^1(\Omega_\ve^i\cup \Omega_\ve^e)$ with the given jump $[\phi_\ve]=w_\ve$ on $\Gamma_\ve$.
Since $V_\ve^0=V_\ve^i - V_\ve^e$ on $\Gamma_\ve$, {by (H1)  we have that}
\begin{align*}
&\int_{\Omega_\ve^i\cup \Omega_\ve^e} a_\ve |\nabla z_\ve|^2\, dx \Big|_{t=0}
= (A_\ve w_\ve, w_\ve)\Big|_{t=0}\\
& \le \int_{\Omega_\ve^i} a_\ve |\nabla V_\ve^i|^2\, dx + \int_{\Omega_\ve^e} a_\ve |\nabla V_\ve^e|^2\, dx\le C.
\end{align*}
The last estimate together with the trace inequality and \eqref{est:extension} completes the proof. 
\end{proof}
Now we prove the a priori estimates for the solution of solution of \eqref{eq:abstract-parabol} for $t\in[0,T]$.
\begin{lemma}[A priori estimates]
\label{lm:apriori-est}
Let $W_\ve=(w_\ve, h_\ve)$ be a solution of \eqref{eq:abstract-parabol}. Then, for $t \in [0,T]$, the following estimates hold:
\begin{enumerate}[(i)]
\item
$\displaystyle \ve \int_{\Gamma_\ve} |w_\ve|^4 \, d\sigma
+ \ve \int_0^t \int_{\Gamma_\ve} |\partial_\tau w_\ve|^2 \, d\sigma \, d\tau \le C$.
\item $\displaystyle \ve \int_{\Gamma_\ve} |h_\ve|^2 \, d\sigma
+ \ve \int_0^t \int_{\Gamma_\ve} |\partial_\tau h_\ve|^2 \, d\sigma \, d\tau\le C$.
\item
Let $z_\ve= e^{-\lambda t} u_\ve$ with the jump $[z_\ve]=w_\ve$ on $\Gamma_\ve$. Then, {one has that}
$$
\displaystyle \int_{\Omega^{i}_{\ve} \cup \Omega^{e}_{\ve}}  (|z_\ve|^2 + |\nabla z_\ve|^2) \, dx \le C,
$$
\end{enumerate}
for a constant $C$ independent of $\ve$, but depending on $T$, and the norms of initial functions $\|G_\ve^0\|_{L^2(\Gamma_\ve)}$, $\|V_\ve^0\|_{L^4(\Gamma_\ve)}$,  $\|V_\ve^l\|_{H^1(\Omega)}$. 
\end{lemma}

\begin{proof}

We {will} work with the equation in 
vector form \eqref{eq:abstract-parabol} and derive the a priori estimates for the pair $(w_\ve, h_\ve)$. Let $z_\ve$ be the solution of the stationary problem with the jump $w_\ve$:
\begin{align}
\label{eq:stat-prob-3}
&-\di{a_\ve \nabla  z_\ve}  =0, \, &x & \in \Omega_\ve^i \cup \Omega_\ve^e, \nonumber\\
&a_e \nabla z_\ve^e\cdot \nu  = a_i \nabla z_\ve^i\cdot \nu, \, &x &\in  \Gamma_\ve, \nonumber\\
&z_\ve^i - z_\ve^e=w_\ve, \, &x& \in  \Gamma_\ve, 
\\
&a_\ve \nabla z_\ve\cdot \nu   = 0, \, &x& \in \Gamma_{m,\ve}, \nonumber \\
&a_e \nabla z_\ve\cdot \nu   = \frac{e^{-\lambda t}}{c_m}J_\ve^e, \, &x& \in \Sigma, \nonumber \\
&z_\ve   = 0, \, &x& \in (S_{0} \cup S_L). \nonumber
\end{align}
We multiply \eqref{eq:abstract-parabol} by $W_\ve$ and integrate over $\Gamma_\ve$:
\begin{align}
&\frac{\ve}{2} \partial_t \int_{\Gamma_\ve} |w_\ve|^2\, d\sigma
+ \frac{1}{c_m} \int_{\Omega_\ve^i\cup \Omega_\ve^e} a_\ve \nabla z_\ve\cdot \nabla z_\ve\, dx
+ \frac{\ve}{c_m} \int_{\Gamma_\ve} \frac{e^{2\lambda t}}{3} w_\ve^4\, d\sigma \nonumber\\
\label{eq:apriori-1}
&+ \ve\left(\lambda-\frac{1}{c_m}\right)\int_{\Gamma_\ve} |w_\ve|^2\, d\sigma
-\ve\left(\theta + \frac{1}{c_m}\right) \int_{\Gamma_\ve} h_\ve w_\ve\, d\sigma
+ \frac{\ve}{2} \partial_t \int_{\Gamma_\ve} |h_\ve|^2\, d\sigma\\
&+ \ve(\lambda+b)\int_{\Gamma_\ve} |h_\ve|^2\, d\sigma
= \frac{e^{-\lambda t}}{c_m}\int_\Sigma J_\ve^e z_\ve\,d\sigma + \ve a e^{-\lambda t} \int_{\Gamma_\ve} h_\ve d\sigma. \nonumber
\end{align}
It is known \cite{acerbi1992extension} that there exists an extension operator $P_\ve$ from $\Omega_\ve^e$ to $\Omega$ such that  $\|\nabla P_\ve z_\ve^e\|_{L^2(\Omega)} \le C \|\nabla z_\ve^e\|_{L^2(\Omega_\ve^e)}$ with a constant $C$ independent of $\ve$.
This result combined with the Friedrichs inequality ($z_\ve=0$ on $S_0\cup S_L$) implies that 
\begin{align}
\label{est:extension}
\|P_\ve z_\ve^e\|_{H^1(\Omega)} \le C\|\nabla z_\ve^e\|_{L^2(\Omega_\ve^e)}. 
\end{align}
By the trace inequality, the $L^2(\Sigma)$-norm of $z_\ve$ is then bounded by $\|\nabla z_\ve^e\|_{L^2(\Omega_\ve^e)}$. 
Using the Young inequality with a parameter in \eqref{eq:apriori-1} and \eqref{est:extension}, yields
\begin{align}
&\partial_t \left(\ve\int_{\Gamma_\ve} |w_\ve|^2 d\sigma
+ \ve \int_{\Gamma_\ve} |h_\ve|^2 d\sigma \right)
+ \int_{\Omega_\ve^i\cup \Omega_\ve^e} |\nabla z_\ve|^2\, dx
+ \ve \int_{\Gamma_\ve} |w_\ve|^4 d\sigma \nonumber\\
&+ \left(\ve\int_{\Gamma_\ve} |w_\ve|^2 d\sigma
+ \ve \int_{\Gamma_\ve} |h_\ve|^2 d\sigma \right)
\le C\int_\Sigma |J_\ve^e|^2\, d\sigma.
\label{eq:apriori-2}
\end{align}
Applying the Gr\"onwall inequality in \eqref{eq:apriori-2}, we obtain the following estimate:
\begin{align}
\label{est:L^2}
\ve \int_{\Gamma_\ve} |w_\ve|^2 d\sigma
+ \ve \int_{\Gamma_\ve} |h_\ve|^2 d\sigma \le C.
\end{align}
Integrating \eqref{eq:apriori-2} with respect to $t$ gives
\begin{align}
\label{est:L^4-grad}
&\int_0^t \int_{\Omega_\ve^i\cup \Omega_\ve^e} |\nabla z_\ve|^2\, dx
+ \ve \int_0^t \int_{\Gamma_\ve} |w_\ve|^4 d\sigma\\
&\le
C\left(\int_0^t \int_{\Sigma}|J_\ve^e|^2 d\sigma d\tau  + \ve \int_{\Gamma_\ve} |V_\ve^0|^2\, d\sigma
+ \ve \int_{\Gamma_\ve} |G_\ve^0|^2\, d\sigma\right).\nonumber
\end{align}
Next,  {we} derive the estimates for $\partial_t W_\ve$. To this end, we multiply \eqref{eq:abstract-parabol} by $\partial_t W_\ve$ and integrate over $(0,t)\times\Gamma_\ve$:
\begin{align}
\nonumber
&\frac{\ve}{2} \int_0^t \int_{\Gamma_\ve} |\partial_\tau w_\ve|^2\, d\sigma d\tau
+ \frac{\ve}{2} \int_0^t \int_{\Gamma_\ve} |\partial_\tau h_\ve|^2\, d\sigma d\tau \nonumber\\
&+ \frac{1}{2 c_m} \int_{\Omega_\ve^i\cup \Omega_\ve^e} a_\ve |\nabla z_\ve|^2\, dx
- \frac{1}{2 c_m} \int_{\Omega_\ve^i\cup \Omega_\ve^e} a_\ve |\nabla z_\ve|^2\, dx \Big|_{t=0} \nonumber\\
&+ \frac{\ve}{12 c_m} e^{2\lambda t}\int_{\Gamma_\ve} |w_\ve|^4\, d\sigma
-\frac{\ve}{12 c_m} \int_{\Gamma_\ve} |V_\ve^0|^4\, d\sigma
\nonumber\\
&+ \frac{\ve}{2}(\lambda-\frac{1}{c_m})\int_{\Gamma_\ve} |w_\ve|^2\, d\sigma
- \frac{\ve}{2}(\lambda-\frac{1}{c_m})\int_{\Gamma_\ve} |V_\ve^0|^2\, d\sigma
\nonumber
\\
\label{eq:apriori-4}
&
+ \frac{\ve}{2}(\lambda+b) \int_{\Gamma_\ve} |h_\ve|^2\, d\sigma
- \frac{\ve}{2} (\lambda+b) \int_{\Gamma_\ve} |G_\ve^0|^2\, d\sigma\\
&
\le 
2\lambda \ve \int_0^t e^{2\lambda \tau} \int_{\Gamma_\ve} |w_\ve|^4\, d\sigma d\tau \nonumber\\
&+ 2\theta^2 \ve \int_0^t \int_{\Gamma_\ve} |w_\ve|^2\, d\sigma d\tau
+ \frac{2\ve}{c_m^2} \int_0^t \int_{\Gamma_\ve} |h_\ve|^2\, d\sigma d\tau \nonumber\\
&+ \frac{e^{-\lambda t}}{c_m}\int_\Sigma J_\ve^e z_\ve\,d\sigma
- \frac{1}{c_m}\int_\Sigma J_\ve^e z_\ve\,d\sigma\Big|_{t=0}\nonumber\\
&+\frac{\lambda}{c_m} \int_0^t e^{-\lambda \tau}\int_\Sigma
J_\ve^e z_\ve\,d\sigma d\tau
- \int_0^t \frac{e^{-\lambda \tau}}{c_m}\int_\Sigma
\partial_\tau J_\ve^e z_\ve\,d\sigma d\tau \nonumber\\
& + \ve a e^{-\lambda t} \int_{\Gamma_\ve} h_\ve \, d\sigma
- \ve a \int_{\Gamma_\ve} G_\ve^0 \, d\sigma
+ \ve a \lambda \int_0^t e^{-\lambda \tau} \int_{\Gamma_\ve} h_\ve \, d\sigma d\tau. \nonumber
\end{align}

Combining \eqref{est:L^2}, \eqref{est:L^4-grad}, and \eqref{est:t=0} we get
\begin{align*}
\ve \int_0^t \int_{\Gamma_\ve} |\partial_\tau w_\ve|^2\, d\sigma d\tau
+ \int_{\Omega_\ve^i\cup \Omega_\ve^e} |\nabla z_\ve|^2\, dx
+\ve \int_{\Gamma_\ve} |w_\ve|^4\, d\sigma \le C.
\end{align*}
Thanks to the homogeneous Dirichlet boundary condition on the bases $S_0 \cup S_L$, the $L^2$-norm of $z_\ve$ is estimated in terms on the $\nabla z_\ve$. Namely,
\begin{align*}
\int_{\Omega_\ve^i} |z_\ve^i|^2\, dx &\le C
\int_{\Omega_\ve^i} |\partial_{x_1} z_\ve^i|^2\, dx,\\
\int_{\Omega_\ve^e} |z_\ve^e|^2\, dx &\le C
\int_{\Omega_\ve^e} |\nabla z_\ve^e|^2\, dx.
\end{align*}
The proof of Lemma \ref{lm:apriori-est} is finally complete.
\end{proof}
\section{Derivation of the 
macroscopic problem}
\label{Rabota_Dlya_Nastoyaschih_Patsanov}
\subsection{Formal asymptotic expansions}
{So as to provide an insight on} how the effective coefficients and the corresponding cell problems in \eqref{eq:hom-prob} appear, we apply {the} formal asymptotic expansion method {to} the stationary problem $A_\ve v_\ve = \ve g$. {Specifically, we write}
\begin{align}
\label{eq:stat-prob-expansion}
-&\di{a_\ve \nabla  u_\ve}  =0, \, &x & \in \Omega_\ve^i \cup \Omega_\ve^e, \nonumber\\
&a_e \nabla u_\ve^e\cdot \nu  = a_i \nabla u_\ve^i\cdot \nu = \ve g(x), \, &x &\in  \Gamma_\ve, \nonumber\\
&u_\ve^i - u_\ve^e=v_\ve, \, &x& \in  \Gamma_\ve, 
\\
&a_e \nabla u_\ve\cdot \nu   = 0, \, &x& \in \Gamma_\ve^m\cup \Sigma, \nonumber \\
& u_\ve   = 0, \, &x& \in (S_{0} \cup S_L). \nonumber
\end{align}
Here $g=g(x)$ is some smooth function. 
Take 
\begin{align*}
u_\varepsilon^l(x) \sim u_0^l(x,y) + \varepsilon u_1^l(x,y) + \varepsilon^2 u_2^l(x,y) + \ldots, \quad y=\frac{x}{\varepsilon},
\end{align*}
where $x \in \Omega_\ve^l$ and $y\in Y_l$, $l \in \{i, e \}$.
Then we get
\begin{align*}
    \mbox{div} (a_l \nabla u_\varepsilon^l ) 
    & \sim \frac{1}{\varepsilon^2}\mbox{div}_y(a_l \nabla_y u_0^l )\\
    &+ \frac{1}{\varepsilon}\left(\mbox{div}_y(a_l \nabla_x u_0^l)  +\mbox{div}_y(a_l \nabla_y u_1^l) +\mbox{div}_x(a_l \nabla_y u_0^l ) \right)\\
    &+ \mbox{div}_x(a_l \nabla_x u_0^l ) 
    + \mbox{div}_x(a_l \nabla_y u_1^l ) 
    + \mbox{div}_y(a_l \nabla_x u_1^l)
    + \mbox{div}_y(a_l \nabla_y u_2^l)\\
     &+ \varepsilon \left( \mbox{div}_x(a_l \nabla_x u_1^l )) + \mbox{div}_x(a_l \nabla_y u_2^l )) + \mbox{div}_y(a_l \nabla_x u_2^l )) \right)\\
     &+ \varepsilon^2 \mbox{div}_x(a_l \nabla_x u_2^l ).
\end{align*}
Taking the terms of order $\varepsilon^{-2}$ in the volume and the ones of order $\varepsilon^{-1}$ on the boundary, we obtain the following problem for $u_0^l$:
\begin{align*}
    - &\mbox{div}_y(a_l \nabla_y u_0^l ) = 0, \quad &y&\in Y_l ,\\
    &a_l \nabla_y u_0^l = 0 \quad &y&\in \Gamma \cup \Gamma^m,\\
    &u_0^i \mbox{ is } 1\mbox{-periodic in }y_1, &&\\
    &\mbox { and } u_0^e \mbox{ is } Y\mbox{-periodic}.&&
\end{align*}
The solution (defined up to an additive constant) does not depend on the fast variable $y$:
\begin{align}\label{eq:asymp-0}
 u_0^l(x,y) = u_0^l(x), \quad l=i, e.   
\end{align}
{For} the next step, we take the terms of order $\varepsilon^{-1}$ in the volume and {those} of order $1$ on the boundary:
\begin{align}
- &\mbox{div}_y(a_l \nabla_y u_1^l) = 0, \quad &y&\in Y_l, \nonumber\\
&a_l \nabla_y u_1^l \cdot \nu = - a_l \nabla_x u_0^l \cdot \nu, \quad &y&\in\Gamma \cup \Gamma_m,
\label{eq:u1}\\ 
&u_1^i \mbox{ is } 1\mbox{-periodic in }y_1 && \nonumber\\
&\mbox { and } u_1^e \mbox{ is } Y\mbox{-periodic}.\nonumber
\end{align}
The solvability condition reads $\displaystyle
    -\int_{\Gamma} a_l \nabla_x u_0^l \cdot \nu = 0$,
which is fulfilled thanks to \eqref{eq:asymp-0}. {By seeking} a solution of \eqref{eq:u1} in the form $u_1^l(x,y) = \mathbf{N}^l(y) \cdot \nabla_x u_0^l (x)$, we obtain
\begin{align*}
a^{l}\nabla_y u_{1}^{l}(x,y) \cdot \nu 
= a^{l} \partial_{y_{j}}N^{l}_{i}(y)\nu_{j}  \partial_{x_i}u_{0}^{l}(x),
\end{align*}
where we assume summation over the repeated indexes. 
The boundary condition in \eqref{eq:u1} yields a boundary condition for $N_i$ on $\Gamma\cup \Gamma_m$:
\begin{align*}
\left( \partial_{y_{j}}N^{l}_{i}(y) + \delta_{i,j}\right)\nu_{j}  = 0.
\end{align*}
\noindent
Then, the functions $N_k^e$, $k=1, 2, 3$, solve the cell problems:
\begin{align}
-&\Delta N_k^e = 0, &y&\in Y_e, \nonumber\\
&\nabla N_k^e \cdot \nu = -\nu_k, \hskip -2.5cm&y&\in \Gamma\cup \Gamma_m,\label{eq:cell-prob-ee}\\
&y \mapsto N_k^e(y) \,\,\,\mbox{is} \,\,\, Y-\mbox{periodic};\hskip -2.5cm&& \nonumber
\end{align}
For the functions $N_k^i$, due to the periodicity in only one variable $y_1$, one can see that $N_k^i(y)=-y_k$ for $k\neq 1$, that yields $\partial_{l\neq k}N_k^i=0$. The first component $N_1^i$ solves the problem 
\begin{align}
-&\Delta N_1^i = 0, &y&\in Y_i, \nonumber\\
&\nabla N_1^i \cdot \nu = -\nu_1, \hskip -2.5cm&y&\in \Gamma\cup \Gamma_m,\label{eq:cell-prob-ii}\\
&y \mapsto N_1^i(y) \,\,\,\mbox{is} \,\,\, 1-\mbox{periodic};\hskip -2.5cm&&\nonumber
\end{align}
Finally, taking the terms of order $1$ in the volume and the ones of order $\epsilon^{1}$ on the boundary, we obtain the following problem for $u_2^{l}$:
\begin{align*}
- &\mbox{div}_{y}(a^{l} \nabla_{y} u_{2}^{l}) = \mbox{div}_{x}(a^{l} \nabla_{x} u_{0}^{l}) + \mbox{div}_{x}(a^{l} \nabla_{y} u_{1}^{l}) + \mbox{div}_{y}(a^{l} \nabla_{x} u_{1}^{l}),\quad &y&\in Y_{l},\\
&a^{l} \nabla_{y} u_{2}^{l} \cdot \nu^l = - a^{l} \nabla_{x} u_{1}^{l} \cdot \nu^l  + g(x), \quad &y&\in \Gamma ,\\ 
& a^{l} \nabla_{y} u_{2}^{l} \cdot \nu=0, \quad &y&\in \Gamma_m,\\
&u_2^i \mbox{ is } 1\mbox{-periodic in }y_1 &&\\
&\mbox { and } u_2^e \mbox{ is } Y\mbox{-periodic}. &&
\end{align*}
Here $\nu^l$ is the exterior unit normal, and $\nu^e = - \nu^i$ on $\Gamma$. The solvability condition {reads}
\begin{align*}
    \int_{Y_{l}} \left( \mbox{div}_{x}(a^{l} \nabla_{x} u_{0}^{l}) + \mbox{div}_{x}(a^{l} \nabla_{y} u_{1}^{l}) + \mbox{div}_{y}(a^{l} \nabla_{x} u_{1}^{l}) \right) dY - \int_{\Gamma}a^{l}\nabla_{x}u_{2}^{l}\cdot \nu^l d\sigma = 0.
\end{align*}
Integrating by parts in the third term of the volume integral, substituting the expression $u_1^l(x,y)= N_i^l(y) \partial_{x_i} u_0^l(x)$, and taking into account that $N_k^i(y)=-y_k$ and $\int_{Y_i} \partial_{l\neq 1}N_1^i dy=0$, we obtain
\begin{align*} 
-\partial_{kj} u_{0}^{e}(x)\int_{Y_e} a^{e} \left(\partial_{j}N_{k}^{e}(y) + \delta_{kj} \right) dy = |\Gamma| g(x),\\
|Y_i| a_i \partial_{11}u_{0}^{i}(x) = |\Gamma| g(x).
\end{align*}
Introducing the effective coefficient
\begin{align*}
    (a_e^\eff)_{kl}= \frac{1}{|\Gamma|}\int_{Y_e} a_e (\partial_{l} N_{k}^{e}(y) + \delta_{kl} ) dy, \quad {k,l=1,2,3,}
\end{align*}
and adding the boundary conditions on $S_0\cup S_L$ and $\Sigma$, we arrive at\\
\begin{align*}
    &\frac{|Y_i|}{|\Gamma|} a_i \partial_{11} u_0^i
    = -a_e^\eff \Delta u_0^e=g(x), \quad &x& \in \Omega, \nonumber\\
   & u_0^{i, e}=0, \quad &x&\in S_0\cup S_L,\\
   & a_e^\eff \nabla u^e\cdot \nu =0, \quad &x& \in \Sigma.\nonumber
\end{align*}

\vspace{2mm}
\subsection{Derivation of the macroscopic problem}
$ $\\
Since the axons inside the bundle are disconnected,
a priori estimates provided by Lemma \ref{lm:apriori-est} do not imply the strong convergence of the transmembrane potential $v_\ve$ on $\Gamma_\ve$. {In turn, this} makes passing to the limit in the nonlinear term $I_{ion}$ problematic. We choose to combine the two-scale convergence machinery with the method of monotone operators due to G. Minty \cite{minty1962}. For reader's convenience we provide a brief description of the method for a simple case in Appendix A, while 
its adaptation for problem \eqref{eq:orig-prob} is presented 
in Section \ref{Rabota_Dlya_Nastoyaschih_Patsanov}.
For passage to the limit, as $\ve \to 0$, we will use the two-scale convergence \cite{Al-1992}. We refer to \cite{AlDa-95} for two-scale convergence on periodic surfaces (namely, on $\Gamma_\ve$).


\begin{definition}
We say that 
a sequence $\{u_\ve^l(t,x)\}$ two-scale converges to the function $u_0^l(t, x, y)$ in $L^2(0,T; L^2(\Omega_{\ve}^l))$, $l=i, e$, as $\ve \to 0$, and write
\begin{equation*}
    u_\ve^l(t,x) \mathop{\rightharpoonup}^2 u_0^l(t, x, y),
\end{equation*}
if
\begin{enumerate}[(i)]
\item
$\displaystyle \int_0^T \int_{\Omega_{\ve}^l} |u_\ve|^2 dx \, dt < C$. 
\item For any $\phi(t,x) \in C(0,T; L^2(\Omega))$, $\psi(y)\in L^2(Y_l)$ we have
\begin{align*}
&\lim\limits_{\ve \to 0} \int_0^T \int_{\Omega_{\ve}^l} u_\ve^l(t, x) \phi(t,x) \psi\left(\frac{x}{\ve}\right)\, dx\, dt\\ 
&=
\frac{1}{|Y|}\int_0^T \int_\Omega \int_{Y^l} u_0^l(t, x,y) \phi(t,x) \psi(y)\, dy \,dx\, dt,
\end{align*}
for some function $u_0^l\in L^2(0,T; L^2(\Omega\times Y))$. 
\end{enumerate} 
\end{definition}

\begin{definition}
A sequence $\{v_\ve(t,x)\}$ converges two-scale to the function $v_0(t, x, y)$ in $L^2(0,T; L^2(\Gamma_\ve))$, as $\ve \to 0$, if
\begin{enumerate}[(i)]
\item
$\displaystyle \ve \int_0^T \int_{\Gamma_\ve} v_\ve^2 \, d\sigma \, dt < C$.
\item For any $\phi(t,x) \in C([0,T]; C(\overline \Omega))$, $\psi(y)\in C(\Gamma)$ we have {that}
\begin{align*}
&\lim\limits_{\ve \to 0}\ve \int_0^T \int_{\Gamma_\ve} v_\ve(t, x) \phi(t,x) \psi\left(\frac{x}{\ve}\right)\, d\sigma_x\, dt\\ 
&=
\frac{1}{|Y|}\int_0^T \int_\Omega \int_{\Gamma} v_0(t, x, y) \phi(t,x) \psi(y)\, d\sigma_y \,dx\, dt
\end{align*}
for some function $v_0 \in L^2(0,T; L^2(\Omega\times \Gamma))$.
\item We say that $\{v_\ve\}$ converges $t$-pointwise two-scale in $L^2(\Gamma_\ve)$ if, for any $t\in[0,T]$, and for any $\phi(x) \in C(\overline \Omega)$, $\psi(y)\in C(\Gamma)$ we have
\begin{align*}
&\lim\limits_{\ve \to 0}\ve \int_{\Gamma_\ve} v_\ve(t, x) \phi(x) \psi\left(\frac{x}{\ve}\right)\, d\sigma_x
=
\frac{1}{|Y|} \int_\Omega \int_{\Gamma} v_0(t, x, y) \phi(x) \psi(y)\, d\sigma_y \,dx
\end{align*}
for some function $v_0 \in L^2(0,T; L^2(\Omega\times \Gamma))$.
\end{enumerate} 
\end{definition}
\vspace{2mm}
\begin{lemma}
\label{lm:compact}
Let $W_\ve$ be a solution of \eqref{eq:abstract-parabol}, and let $z_\ve$ be a solution of problem \eqref{eq:stat-prob-3}. Then there exist functions $z_0^l \in L^2(0,T; L^2(\Omega))$, $l=i, e$, such that $\partial_{x_1} z_0^i, \partial_{x_j} z_0^e \in L^2(0,T; L^2(\Omega))$ ($j=1, 2, 3$), $w_0=z_0^i-z_0^e \in L^4(0,T; L^4(\Omega))$, and up to a subsequence, as $\ve \to 0$, the following two-scale convergence holds:
\begin{itemize}
\item[(i)] $\displaystyle\chi^l\left(\frac{x}{\ve}\right)z_\ve^l(t,x) \quad 
\overset{\,2\,}{\rightharpoonup} \quad \chi^l(y) z_0^l(t, x)$ in $L^2(0,T; L^2(\Omega_\ve^l))$, $l=i, e$.\\[2mm]
\item[(ii)] $\displaystyle\chi^i\left(\frac{x}{\ve}\right) \nabla z_\ve^i(t,x) \quad \mathop{\rightharpoonup}^2 \quad
\chi^i(y) \big[{\bf e}_1 \partial_{x_1} z_0^i(t, x) + \nabla_y z_1^i(t, x,y)\big]$, where $z_1^i(t, x,y) \in L^2((0,T)\times \Omega; H^1(Y_i))$ is $1$-periodic in $y_1$.\\[2mm]
\item[(iii)] $\displaystyle\chi^e\left(\frac{x}{\ve}\right) \nabla z_\ve^e(t,x) \quad \mathop{\rightharpoonup}^2 \quad
\chi^e(y) \big[\nabla z_0^e(t,x) + \nabla_y z_1^e(t, x,y)\big]$, where $z_1^e(t, x,y) \in L^2((0,T)\times \Omega; H^1(Y_e))$ is $Y$-periodic in $y$.\\[2mm]
\item[(iv)] $\displaystyle w_\ve \,\, \mathop{\rightharpoonup}^2 \,\, w_0(t,x)$ $t$--pointwise in $L^2(\Gamma_\ve)$, and $w_0=(z_0^i - z_0^e)$. Moreover, $\displaystyle \partial_t w_\ve \,\, \mathop{\rightharpoonup}^2 \,\, \partial_t w_0$ in $L^2(0,T; L^2(\Gamma_\ve))$. \\[2mm]
\item[(v)] $\displaystyle h_\ve \,\, \mathop{\rightharpoonup}^2 \,\, \widetilde h_0(t, x,y)$ $t$--pointwise in $L^2(\Gamma_\ve)$, and $\displaystyle \partial_t h_\ve \,\, \mathop{\rightharpoonup}^2 \,\, \partial_t \widetilde h_0$ in $L^2(0,T; L^2(\Gamma_\ve))$.
\end{itemize}
\end{lemma}

\begin{proof}
From a priori estimates the two-scale convergence of $z_\ve^e$ and $\nabla z_\ve^e$ is proved applying standard arguments (see \cite{Al-1992}). When it comes to $z_\ve^i$ and its gradient, the main difficulty stems
from the fact that $\Omega_\ve^i$ consists of many disconnected components. 

Since $z_\ve^i$ is bounded uniformly in $\ve$ {(cf.~ Lemma \ref{lm:apriori-est})} in $L^2((0,T)\times \Omega_\ve^i)$, there exists a subsequence---still denoted by $\{z_\ve^i\}$---such that $\chi^i(\frac{x}{\ve})z_\ve^i(t, x)$ converging two-scale to some $\chi^i(y)z_0^i(t, x,y)$ in $L^2(0,T; L^2(\Omega\times Y))$. Similarly, due to \eqref{est:L^4-grad}, up to a subsequence, $\chi^i\left(\frac{x}{\ve}\right) \nabla z_\ve^i(t, x)$ converges two-scale to $\chi^i(y) p^i(t, x,y)$. Let us show that $z_0^i=z_0^i(t, x)$. Take a smooth test function $\displaystyle\Phi\left(t, x, \frac{x}{\ve}\right)= \varphi(t, x) \psi\left(\frac{x}{\ve}\right)$, where $\varphi\in C([0,T]; C_0^\infty(\Omega))$, and $\psi \in (C^\infty(Y_i))^3$ is $1$-periodic in $y_1$ and such that $\psi=0$ on $\Gamma_{mi}\cup \Gamma$.
\begin{align*}
\ve \int_0^T\int_{\Omega_\ve^i} \nabla z_\ve^i(t, x)\cdot \varphi(t, x) \psi\left(\frac{x}{\ve}\right)\, dx dt
= &-\ve \int_0^T \int_{\Omega_\ve^i} z_\ve^i(t,x)\nabla \varphi(t,x) \cdot \psi\left(\frac{x}{\ve}\right)\, dx dt\\
&- \int_0^T \int_{\Omega_\ve^i} z_\ve^i(t,x) \varphi(t,x){\rm div}_y \psi\left(\frac{x}{\ve}\right)\, dx dt.
\end{align*}
Passing to the limit, {we derive}
\begin{align*}
\frac{1}{|Y|} \int_0^T \int_{\Omega} \int_{Y_i} z_0^i(t,x,y) \varphi(t,x) {\rm div}_y \psi(y)\, dy dx dt=0,
\end{align*}
{which} implies {that} $\partial_{y_i} z_0^i(t,x,y)=0$, $i=1, 2, 3$. Thus, $z_0^i=z_0^i(t,x)$.

Next we prove that $\partial_{x_1} z_0^i \in L^2((0,T) \times \Omega)$. Let us take a test function $\Phi\left(t, x, \frac{x}{\ve}\right)= \varphi(t, x) {\mathbf e}_1 + \varphi(t, x) \nabla_y N_1^i\left(\frac{x}{\ve}\right)$ such that
\begin{align}
    &\Delta_y N_1^i=0, \quad Y_i, \nonumber\\
    &\nabla N_1^i\cdot \nu = -\nu_1, \quad \Gamma\cup \Gamma_{mi},
    \label{eq:N1}\\
    &N_1^i \,\, \mbox{is 1-periodic in}\,\, y_1. \nonumber
\end{align}
Integrating by parts yields
\begin{align*}
&\int_0^T \int_{\Omega_\ve^i} \nabla z_\ve^i(t,x) \cdot \Phi\left(t, x, \frac{x}{\ve}\right)\, dx dt\\
&= -\int_0^T \int_{\Omega_\ve^i} z_\ve^i(t,x)\left({\mathbf e}_1 + \nabla_y N_1^i\left(\frac{x}{\ve}\right)\right) \cdot \nabla \varphi(t,x)\, dx dt,
\end{align*}
and  passing to the limit, as $\ve \to 0$, we obtain
\begin{align}
\label{eq:nomerok}
&\frac{1}{|Y|} \int_0^T \int_{\Omega} \int_{Y_i} p^i(t,x,y) \cdot \varphi(t,x)\left({\mathbf e}_1 + \nabla_y N_1^i(y)\right)\, dy dx dt\\
\nonumber
&= -\frac{1}{|Y|} \int_0^T\int_{\Omega} \int_{Y_i}
z_0^i(t,x) \nabla \varphi(t,x) \cdot \left({\mathbf e}_1 + \nabla_y N_1^i(y)\right)\, dy dx dt.
\end{align}
Let us observe that $\int_{Y_i} \partial_{y_k} N_1^i(y) \, dy=0$ for $k\neq 1$. Indeed, for $k\neq 1$, $y_k$ can be taken as a test function in \eqref{eq:N1}:
\begin{align*}
0=-\int_{Y_i} \Delta N_1^i(y) y_k\, dy
= \int_{Y_i} \partial_{y_k} N_1^i(y) \, dy.
\end{align*}
{Furthermore, it holds that} $$\int_{Y_i} \left(\delta_{1k}+ \partial_{y_k} N_1^i(y)\right)\, dy= \delta_{1k}|\Gamma|\frac{a_i^\eff}{a_i}.$$ 
{Consequently,} it is straightforward to check that 
\begin{align}
\label{eq:a_i^eff}
a_i^\eff = \frac{1}{|\Gamma|}\int_{Y_i} a_i\left(1+ \partial_{y_1} N_1^i(y)\right)\, dy=
\frac{1}{|\Gamma|} \int_{Y_i} a_i\left(1+ \partial_{y_1} N_1^i(y)\right)^2\, dy >0.
\end{align}
We turn back to \eqref{eq:nomerok}. Due to \eqref{eq:a_i^eff}, {we have the estimate} 
\begin{align*}
&\left|\int_0^T \int_{\Omega} z_0^i(t,x) \partial_{x_1} \varphi(t,x)\, dx dt\right|\\
&= \left| \frac{a_i}{(a_i^\eff)_{11}} 
\int_0^T \int_{\Omega} \int_{Y_i} p^i(t,x,y) \cdot \varphi(t, x)\left({\mathbf e}_1 + \nabla_y N_1^i(y)\right)\, dy dx dt\right|\\[2mm]
&\le C \|\varphi\|_{L^2((0,T)\times\Omega)}.
\end{align*}
Next, we show that $p^i(t, x,y)={\mathbf e}_1 \partial_{x_1} z_0^i(t,x) + \nabla_y z_1^i(t,x,y)$ for some $z_1^i$ periodic in $y_1$. Take a smooth test function $\varphi(t,x)\psi(y)$ such that ${\rm div}_y \psi=0$ in $Y_i$, $\psi\cdot \nu=0$ on $\Gamma_{mi} \cup \Gamma$, and periodic in $y_1$.
\begin{align*}
\int_0^T \int_{\Omega_\ve^i} \nabla z_\ve^i \cdot \varphi(t,x) \psi\left(\frac{x}{\ve}\right) \, dx dt
= - \int_0^T \int_{\Omega_\ve^i} z_\ve^i \nabla \varphi(t,x) \cdot \psi\left(\frac{x}{\ve}\right) \, dx dt.
\end{align*}
Passing to the limit, as $\ve \to 0$ we obtain
\begin{align*}
\frac{1}{|Y|} \int_0^T \int_\Omega \int_{Y_i} p^i \cdot \varphi(t,x) \psi(y) \, dy dx dt
=- \frac{1}{|Y|} \int_0^T \int_\Omega \int_{Y_i} 
z_0^i \nabla \varphi(t,x) \cdot \psi(y)  \, dy dx dt.
\end{align*}
Since $\int_{Y_i} \psi_k(y)\, dy=0$ for $k\neq 1$, 
\begin{align*}
\int_0^T \int_\Omega \int_{Y_i} p^i(t,x,y) \cdot \varphi(t,x) \psi(y) \, dy dx dt
= \int_0^T \int_\Omega \int_{Y_i} 
\partial_{x_1} z_0^i(t,x) \varphi(t,x) \psi_1(y)  \, dy dx dt,
\end{align*}
and thus
\begin{align*}
\int_0^T \int_\Omega \int_{Y_i}
\left(p^i(x,y) - {\mathbf e}_1 \partial_{x_1} z_0^i(t,x)\right) \varphi(t,x) \cdot \psi(y)\, dy dx dt=0.
\end{align*}
Since $\psi$ is solenoidal, there exists $z_1^i(t,x,y) \in L^2((0,T)\times \Omega; H^1(Y_i))$, $1$-periodic in $y_1$, such that
\begin{align*}
p^i(t,x,y)= {\mathbf e}_1 \partial_{x_1} z_0^i(t,x)
+ \nabla_y z_1^i(t,x,y).
\end{align*}

Next we prove that the jump $w_\ve$ converges two-scale in $L^2(0,T; L^2(\Gamma_\ve))$ to $z_0^i- z_0^e$. To this end, for $\psi \in H^{1/2}(\Gamma)$, we consider test functions $\widetilde \psi^l$, $l=i, e$, solving
\begin{align*}
&\Delta \widetilde \psi^l = \frac{1}{|Y_l|} \int_\Gamma \psi\, d\sigma, \quad y\in Y_l, \nonumber\\
&\nabla \widetilde \psi^l \cdot \nu^l = \psi, \quad y\in \Gamma;
\quad \nabla \widetilde \psi^l \cdot \nu^l = 0, \quad y\in \Gamma_{ml},\\
\nonumber
&\widetilde \psi^l \,\, \mbox{is}\,\, Y-\mbox{periodic}.
\end{align*}
Integration by parts yields
\begin{align*}
&\ve \int_0^T\int_{\Gamma_\ve} w_\ve\, \varphi(t,x) \psi\left(\frac{x}{\ve}\right)\, dx dt\\
&= 
\ve \int_0^T \int_{\Omega_\ve^i} \nabla z_\ve^i \cdot \varphi(t,x) \nabla_y \widetilde \psi^i\left(\frac{x}{\ve}\right)\, dx dt
+ \ve \int_0^T\int_{\Omega_\ve^i} z_\ve^i\, \nabla \varphi(t,x) \cdot \nabla_y \widetilde \psi^i\left(\frac{x}{\ve}\right) \, dx dt\\
&+ \frac{1}{|Y_i|}\int_0^T \int_{\Omega_\ve^i} z_\ve^i \varphi(t,x) \int_\Gamma \psi(y) \, d\sigma dx dt\\
&
-
\ve \int_0^T\int_{\Omega_\ve^e} \nabla z_\ve^e \cdot \varphi(t,x) \nabla_y \widetilde \psi^e\left(\frac{x}{\ve}\right)\, dx dt
- 
\ve \int_0^T\int_{\Omega_\ve^e} z_\ve^e \nabla \varphi(t,x) \cdot \nabla_y \widetilde \psi^e\left(\frac{x}{\ve}\right) \, dx dt\\
&- \frac{1}{|Y_e|}\int_0^T \int_{\Omega_\ve^e} z_\ve^e \varphi(t,x) \int_\Gamma \psi(y) \, d\sigma dx dt.
\end{align*}
Passing to the limit, as $\ve \to 0$, we get
\begin{align*}
&\frac{1}{|Y|}\int_0^T \int_{\Omega} \int_\Gamma w_0(t, x,y) \varphi(t, x) \psi(y) \, d\sigma dx dt\\
&= \frac{1}{|Y|} \int_0^T \int_\Omega \int_\Gamma
(z_0^i- z_0^e) \varphi(t,x) \psi(y) \, d\sigma dx dt,
\end{align*}
that proves the two-scale convergence of $w_\ve$ to the difference $w_0=z_0^i-z_0^e$. 

Note that the uniform bound of $w_\ve$ in $L^4((0,T)\times \Gamma_\ve)$---by Lemma \ref{lm:apriori-est}(i)---implies $w_0\in L^4((0,T)\times \Omega)$. Indeed, for smooth $\varphi(t, x)$, we have {that}
\begin{align*}
&|\Gamma|\int_0^T \int_\Omega
w_0(t, x) \varphi(t, x)\, dx dt 
= \lim \limits_{\ve \to 0}
\ve |Y| \int_0^T \int_{\Gamma_\ve} 
w_\ve(t, x) \varphi(t, x)\, d\sigma dt\\
&\le |Y| \lim \limits_{\ve \to 0} \left(\ve \int_0^T \int_{\Gamma_\ve} |w_\ve|^4\, d\sigma dt\right)^\frac{1}{4} 
\left(\ve \int_0^T \int_{\Gamma_\ve} |\varphi(t, x)|^{4/3}\, d\sigma dt\right)^\frac{3}{4}\\
& \le C \lim \limits_{\ve \to 0} \left(\ve \int_0^T \int_{\Gamma_\ve} |\varphi(t, x)|^{\frac{4}{3}}\, d\sigma_{x} dt\right)^\frac{3}{4} \\
&= C \left(\frac{|\Gamma|}{|Y|}\int_0^T \int_{\Omega} \int_\Gamma |\varphi(t, x)|^{\frac{4}{3}}\, dx dt\right)^\frac{3}{4}.
\end{align*}
By density of smooth functions in $L^\frac{4}{3}((0,T)\times \Omega)$, $\|w_0\|_{L^4((0,T)\times \Omega)} \le C$. \\[2mm]
Thanks to the uniform in $\ve$ estimate \eqref{est:L^4-grad}, $(v)$ and $(vi)$ hold. 
\end{proof}
\begin{lemma}
\label{lm:initial-condition}
Let the initial functions $V_\ve^0$ satisfy hypothesis (H1). Then, {one has that}
\begin{align*}
\ve \int_{\Gamma_\ve} |V_\ve^0|^2\, d\sigma= 
\frac{|\Gamma|}{|Y|} \int_{\Omega} |V^i - V^e|^2 \, dx.
\end{align*}
\end{lemma}
\begin{proof}
Approximating $V^0= V^i-V^e$ by smooth functions $V_\delta^0$ in $H^1(\Omega)$, we find
\begin{align}
&\ve \int_{\Gamma_\ve} |V_\ve^0|^2\, d\sigma
= 
\ve \int_{\Gamma_\ve} |V_\ve^0 - V_\delta^0|^2\, d\sigma \nonumber\\
\label{aux-2}
&+ 2 \ve \int_{\Gamma_\ve} (V_\ve^0-V_\delta^0) V_\delta^0\, d\sigma
+ \ve \int_{\Gamma_\ve} |V_\delta^0|^2\, d\sigma.
\end{align}
Applying the trace inequality in the rescaled periodicty cell $\ve Y$, adding up over all the cells in $\Omega$, and using assumption (H1) {leads to}
\begin{align*}
&\ve \int_{\Gamma_\ve} |V_\ve^0 - V_\delta^0|^2\, d\sigma
\le 
C\ve^2  \int_\Omega |\nabla (V_\ve^0 - V_\delta^0)|^2\, dx
+ C \int_\Omega |V_\ve^0 - V_\delta^0|^2\, dx\\
&\le
C\ve^2  \int_\Omega |\nabla (V_\ve^0 - V_\delta^0)|^2\, dx
+ C \int_\Omega |V_\ve^0 - V^0|^2\, dx\\
&+ C \int_\Omega |V_\delta^0 - V^0|^2\, dx \quad \to \quad 0, \quad \ve, \delta \to 0.
\end{align*}
Then, since $V_\delta^0$ is smooth, it converges strongly two-scale, and passing to the limit as $\ve \to 0$ in \eqref{aux-2} we obtain
\begin{align*}
\lim_{\delta \to 0}\limsup_{\ve \to 0} \ve \int_{\Gamma_\ve} |V_\ve^0|^2\, d\sigma =
\frac{|\Gamma|}{|Y|} \int_\Omega |V^0|^2\, dx,
\end{align*}
{as stated}.
\end{proof}
\vspace{3mm}

We proceed with the Minty method for passing to the limit in the microscopic problem. Consider arbitrary functions $\mu_0^l(t,x) \in C^\infty([0,T]\times \overline{\Omega})$ and $\mu_1^l(t,x,y)\in C^\infty([0,T]\times \overline{\Omega} \times Y)$, $Y$-periodic in $y$, and such that $\mu_0^l = \mu_1^l = 0$ when $x\in S_0\cap S_L$.
Take the test function 
\begin{align*}
M_\ve &:=
\begin{pmatrix}
[\mu_\ve] \\ \rho
\end{pmatrix}, \quad \mbox{where}\,\, \rho=\rho(t,x),\ \text{and} \\
\mu_\ve(x)&:=
\begin{cases}
\displaystyle
\mu_0^e(t,x)+\ve \mu_{1}^e\left(t,x,\frac{x}{\ve}\right), \quad x\in \Omega_\ve^e \\[2mm]
\displaystyle
\mu_0^i(t,x)+\ve \mu_{1}^i\left(t,x,\frac{x}{\ve}\right),\quad x\in \Omega_\ve^e.
\end{cases}
\end{align*}
{The} monotonicity property of the operator $\mathbb A_\ve(t, \cdot)$ {entails}
\begin{align*}
\int_0^t \int_{\Gamma_\ve}\left(\mathbb A_\ve(\tau, W_\ve) - \mathbb A_\ve(\tau, M_\ve)\right)\cdot\left(W_\ve - M_\ve\right)\, d\sigma d\tau \ge 0.
\end{align*}
It follows from \eqref{eq:abstract-parabol} and the definition of the operator $\mathbb A_\ve(t, \cdot)$ that
\begin{align}
&\ve \int_0^t \int_{\Gamma_\ve}\partial_\tau w_\ve ( [\mu_\ve]-w_\ve)\, d\sigma d\tau
+
\ve \int_0^t\int_{\Gamma_\ve} \partial_\tau h_\ve (\rho-h_\ve)\, d\sigma d\tau \nonumber\\
&+ \frac{1}{c_m}\int_0^t  \int_{\Gamma_\ve} A_\ve[\mu_\ve]([\mu_\ve] - w_\ve)\, d\sigma d\tau
+ 
\ve(\lambda - \frac{1}{c_m})\int_0^t\int_{\Gamma_\ve} [\mu_\ve] ([\mu_\ve]-w_\ve)\, d\sigma d\tau \nonumber\\
\label{reznik-1}
&-\frac{\ve}{c_m} \int_0^t \int_{\Gamma_\ve} \rho ([\mu_\ve]-w_\ve)\, d\sigma d\tau
+\ve(b+\lambda) \int_0^t\int_{\Gamma_\ve} \rho (\rho - h_\ve)\, d\sigma d\tau
\\
\nonumber
&-\ve \theta \int_0^t \int_{\Gamma_\ve} [\mu_\ve] (\rho - h_\ve)\, d\sigma d\tau
+ \ve \frac{1}{3c_m} \int_0^t e^{2\lambda \tau}\int_{\Gamma_\ve}[\mu_\ve]^3 ([\mu_\ve]- w_\ve)\, d\sigma d\tau\\
\nonumber
&+ \int_0^t \int_{\Gamma_\ve} \frac{e^{-\lambda \tau}}{c_m} (a_i \nabla q_\ve \cdot \nu) ([\mu_\ve]-w_\ve)\, d\sigma d\tau
-\ve a\int_0^t \int_{\Gamma_\ve} e^{-\lambda \tau} (\rho - h_\ve)\, d\sigma d\tau \ge 0.
\end{align}
Consider the first two terms
in \eqref{reznik-1}, 
specifically integrals $\ve \int_0^t \int_{\Gamma_\ve}w_\ve \partial_\tau w_\ve \, d\sigma d\tau$
and 
$\ve\int_0^t\int_{\Gamma_\ve} h_\ve \partial_\tau h_\ve\, d\sigma d\tau$.
Integrating by parts with respect to time, passing to the limit as $\ve \to 0$, and using the lower semi-continuity of 
$L^2$-norm with respect to two-scale convergence (Proposition 2.5, \cite{AlDa-95}) and Lemma \ref{lm:initial-condition} {renders}
\begin{align*}
&\limsup \limits_{\ve \to 0}
\Big[\ve \int_0^t \int_{\Gamma_\ve} w_\ve \partial_\tau w_\ve \, d\sigma d\tau
- \frac{|\Gamma|}{|Y|} \int_0^t \int_\Omega w_0 \partial_\tau w_0\, dx d\tau\Big]\\
&
= \limsup \limits_{\ve \to 0}
\Big[\frac{\ve}{2} \int_{\Gamma_\ve} w_\ve^2\, d\sigma \Big|_{\tau=t} - 
\frac{|\Gamma|}{2|Y|} \int_\Omega w_0^2\, dx\Big]\\
&+
\lim \limits_{\ve \to 0}
\Big[
-\frac{\ve}{2} \int_{\Gamma_\ve} (V_\ve^0)^2\, d\sigma
+\frac{|\Gamma|}{2|Y|} \int_{\Omega} (V^0)^2\, dx
\Big]\ge 0.
\end{align*}
Similarly, for the integral of $h_\ve \partial_\tau h_\ve$, denoting the mean value of the two-scale limit $\widetilde h_0(t,x,y)$ in $y$ by $h_0(t,x)=\frac{1}{|\Gamma|}\int_\Gamma \widetilde h_0(t,x,y)\, dy$, we get
\begin{align*}
&\limsup \limits_{\ve \to 0}
\Big[\ve \int_0^t \int_{\Gamma_\ve} h_\ve \partial_\tau h_\ve \, d\sigma d\tau
- \frac{|\Gamma|}{|Y|} \int_0^t \int_\Omega h_0 \partial_\tau h_0\, dx d\tau\Big]\\
&
= \limsup \limits_{\ve \to 0}
\Big[\frac{\ve}{2} \int_{\Gamma_\ve} h_\ve^2\, d\sigma \Big|_{\tau=t} - 
\frac{|\Gamma|}{2|Y|} \int_\Omega h_0^2\, dx \Big|_{\tau=t}\Big]\\
&+
\lim \limits_{\ve \to 0}
\Big[
-\frac{\ve}{2} \int_{\Gamma_\ve} (G_\ve^0)^2\, d\sigma
+\frac{|\Gamma|}{2|Y|} \int_{\Omega} (G^0)^2\, dx
\Big]\ge 0.
\end{align*}
For smooth $\mu_0^l(t,x)$ and  $\mu_{1}^l(t,x,y)$, $l=i, e$, we use Lemma \ref{lm:compact} 
to pass to the limit in the third term:
\begin{align*}
&\frac{1}{c_m}\int_0^t \int_{\Gamma_\ve} A_\ve[\mu_\ve]([\mu_\ve] - w_\ve)\, d\sigma d\tau \nonumber\\
&\to
\frac{1}{c_m|Y|}\int_0^t\int_\Omega \int_{Y_i}
a_i (\nabla \mu_0^i+\nabla_y \mu_{1}^i)\cdot
( \nabla \mu_0^i+\nabla_y \mu_{1}^i-\partial_1 z_0^i{\bf e}_1-\nabla_y z_{1}^i) dx dy d\tau
\\
&+ \frac{1}{c_m|Y|}\int_0^t\int_\Omega \int_{Y_e} 
a_e (\nabla \mu_0^e+\nabla_y \mu_{1}^e)\cdot
( \nabla \mu_0^e+\nabla_y \mu_{1}^e-\nabla z_0^e-\nabla_y z_{1}^e) dx dy d\tau.\nonumber
\end{align*}
Taking the limit in \eqref{reznik-1} 
as $\ve\to 0$ (along a subsequence) we obtain
\begin{align}
&\limsup \limits_{\ve \to 0}
\Big[\ve \int_0^t \int_{\Gamma_\ve} w_\ve \partial_\tau w_\ve \, d\sigma d\tau
- \frac{|\Gamma|}{|Y|} \int_0^t \int_\Omega w_0 \partial_\tau w_0\, dx d\tau\Big] \nonumber\\
&+
\limsup \limits_{\ve \to 0}
\Big[\ve \int_0^t \int_{\Gamma_\ve} h_\ve \partial_\tau h_\ve \, d\sigma d\tau
- \frac{|\Gamma|}{|Y|} \int_0^t \int_\Omega h_0 \partial_\tau h_0\, dx d\tau\Big] \nonumber\\
&\le \frac{|\Gamma|}{|Y|}\int_0^t \int_\Omega\partial_\tau w_0 ([\mu_0]-w_0)\, dx d\tau
+
\frac{|\Gamma|}{|Y|}\int_0^t\int_\Omega\partial_\tau h_0 (\rho-h_0)\, dx d\tau \nonumber\\
&
+\frac{1}{c_m|Y|}\int_0^t\int_\Omega \int_{Y_i}
a_i (\nabla \mu_0^i+\nabla_y \mu_{1}^i)\cdot
( \nabla \mu_0^i+\nabla_y \mu_{1}^i-\partial_1 z_0^i{\bf e}_1-\nabla_y z_{1}^i) dx dy d\tau \nonumber
\\
&+ \frac{1}{c_m|Y|}\int_0^t\int_\Omega \int_{Y_e} 
a_e (\nabla \mu_0^e+\nabla_y \mu_{1}^e)\cdot
( \nabla \mu_0^e+\nabla_y \mu_{1}^e-\nabla z_0^e-\nabla_y z_{1}^e) dxdy d\tau \nonumber\\
\label{aux-1}
&+ (\lambda - \frac{1}{c_m}) \frac{|\Gamma|}{|Y|}\int_0^t \int_\Omega [\mu_0]([\mu_0]-w_0)\, dx d\tau\\
&-\frac{|\Gamma|}{|Y|c_m}
\int_0^t \int_\Omega \rho ( [\mu_0]-w_0)\, dx d\tau
+(b+\lambda)\frac{|\Gamma|}{|Y|} 
\int_0^t \int_\Omega \rho (\rho - h_0)\, dx d\tau
\nonumber \\
&- \theta\frac{|\Gamma|}{|Y|} \int_0^t\int_\Omega [\mu_0] (\rho - h_0)\, dx d\tau
+\frac{1|\Gamma|}{3c_m |Y|} \int_0^t\int_\Omega e^{2\lambda \tau} [\mu_0]^3([\mu_0]- w_0)\, dx d\tau \nonumber\\
&-\int_0^t \int_\Sigma \frac{e^{-\lambda \tau}}{c_m} J^e (\mu_0^e-z_0^e)\, d\sigma d\tau
- a \frac{|\Gamma|}{|Y|} 
\int_0^t \int_{\Omega} e^{-\lambda \tau} (\rho - h_0)\, d\sigma d\tau,\nonumber
\end{align}
where $[\mu_0]=\mu_0^i-\mu_0^e$. 
Consider the spaces
\begin{align*}
H_i& = \{z^i\in L^2(\Omega): \,\,\partial_{x_1} z^i \in L^2(\Omega), \,\, z^i=0 \,\, \mbox{on} \,\, S_0\cup S_L\},\\
H_e& = \{z^e\in L^2(\Omega): \,\,  \nabla z^e \in L^2(\Omega)^3, \,\, z^e=0 \,\, \mbox{on} \,\, S_0\cup S_L\},\nonumber
\end{align*}
with the standard $H^1$-norm in $H_e$, and 
\begin{align*}
\|z\|_{H_i} = \left(\int_\Omega |z|^4 \, dx\right)^\frac{1}{4}
+ \left(\int_\Omega |\partial_{x_1} z|^2 \, dx\right)^\frac{1}{2}.
\end{align*}
By density of smooth functions, inequality \eqref{aux-1} still holds for test functions $\mu_1^l \in L^2((0,T)\times \Omega; H^1(Y_l))$, and  $\mu_0^l \in L^2(0,T; H_l)$ such that $[\mu_0] \in L^4((0,T) \times \Omega)$.  

Modifying the test function $\mu_{1}^i$ by setting 
$\mu_{1}^i(x,y)=\widetilde \mu_{1}^i(x,y)-\nabla_{x^\prime}\mu_0^i\cdot y^\prime$ 
we transform the integrand  in second line of
\eqref{aux-1} to the form 
$$
a_i (\partial_{x_1} \mu_0^i {\bf e}_1+\nabla_y \widetilde \mu_{1}^i)\cdot
( \partial_{x_1} \mu_0^i {\bf e}_1+\nabla_y \widetilde \mu_{1}^i-\partial_{x_1} z_0^i{\bf e}_1-\nabla_y z_{1}^i).
$$
Then, for smooth test 
functions $\psi^{l}(t, x)$, $\varphi(t,x)$ vanishing at $x=0, L$, and
$\Psi^{l}(t,x,y)$ periodic in $y$ and equal to zero when $x=0, L$, $l=i, e$,  we can 
set
\begin{align*}
  \mu_0^{l}(t, x)&=z_0^{l}(t, x)+\delta \psi^{l}(t, x), \quad l=i,e,\\
  \mu_{1}^{e}(t,x,y)&=z_{1}^{e}(t,x,y)+\delta \Psi^{e}(t, x,y),\\
\widetilde \mu_{1}^{i}(t, x,y)&=z_{1}^{i}(t,x,y)+\delta \Psi^{i}(t, x,y),\\
\rho(t,x)&= h_0(t,x) + \delta \varphi(t,x),
\end{align*}
where  $\delta$ is a small auxiliary parameter. Setting $[\psi]=\psi^i-\psi^e$, we have {that}
\begin{align}
&\limsup \limits_{\ve \to 0}
\Big[\ve \int_0^t \int_{\Gamma_\ve} w_\ve \partial_\tau w_\ve \, d\sigma d\tau
- \frac{|\Gamma|}{|Y|} \int_0^t \int_\Omega w_0 \partial_\tau w_0\, dx d\tau\Big] \nonumber\\
&+
\limsup \limits_{\ve \to 0}
\Big[\ve \int_0^t \int_{\Gamma_\ve} h_\ve \partial_\tau h_\ve \, d\sigma d\tau
- \frac{|\Gamma|}{|Y|} \int_0^t \int_\Omega h_0 \partial_\tau h_0\, dx d\tau\Big] \nonumber\\
&\le \frac{\delta |\Gamma|}{|Y|}\int_0^t \int_\Omega\partial_\tau w_0 [\psi]\, dx d\tau
+
\frac{\delta|\Gamma|}{|Y|}\int_0^t\int_\Omega\partial_\tau h_0 \, \varphi\, dx d\tau \nonumber\\
&
+\frac{\delta }{c_m|Y|}\int_0^t\int_\Omega \int_{Y_i} 
a_i (\partial_{x_1} (z_0^i+\delta \psi^{i}) {\bf e}_1+\nabla_y (z_{1}^i +\delta  \Psi^{i}))\cdot
( \partial_{x_1} \psi^{i}{\bf e}_1 + \nabla_y \Psi^{i}) dx dy  d\tau  \nonumber
\\
&+ \frac{\delta}{c_m|Y|}\int_0^t\int_\Omega \int_{Y_e} 
a_e (\nabla (z_0^e+\delta \psi^{e})+\nabla_y (z_{1}^e+\delta \Psi^{e}))\cdot
( \nabla \psi^{e}+\nabla_y \Psi^{e}) dx dy d\tau  \nonumber\\
\label{eq:osnovopolagaushaya}
&+ (\lambda - \frac{1}{c_m}) \frac{\delta|\Gamma|}{|Y|}\int_0^t \int_\Omega (w_0 + \delta [\psi])[\psi]\, dx d\tau\\
&-\frac{\delta|\Gamma|}{|Y|c_m}
\int_0^t \int_\Omega (h_0 + \delta \varphi) [\psi]\, dx d\tau
+(b+\lambda)\frac{\delta|\Gamma|}{|Y|} 
\int_0^t \int_\Omega (h_0 + \delta \varphi) \varphi\, dx d\tau  \nonumber
\\
&- \theta\frac{\delta|\Gamma|}{|Y|} \int_0^t\int_\Omega (w_0 + \delta [\psi]) \varphi\, dx d\tau
+\frac{|\Gamma|}{3c_m |Y|} \delta \int_0^t\int_\Omega e^{2\lambda \tau}(w_0+ \delta[\psi])^3 [\psi]\, dx d\tau \nonumber\\
&-\frac{\delta}{c_m}  \int_0^t \int_\Sigma e^{-\lambda \tau}J^e \psi^e\, d\sigma d\tau
-  a \frac{\delta |\Gamma|}{|Y|} 
\int_0^t \int_{\Omega} e^{-\lambda \tau} \varphi\, d\sigma d\tau. \nonumber
\end{align}
Since the left-hand side of \eqref{eq:osnovopolagaushaya} is non-negative and $\delta$ is arbitrary, we obtain
\begin{align*}
&\limsup_{\ve \to 0}
\Big[\ve \int_{\Gamma_\ve} |w_\ve|^2\, d\sigma
- \frac{|\Gamma|}{|Y|} \int_\Omega |w_0|^2\, dx\Big]=0,\\
&\limsup_{\ve \to 0} \Big[\ve \int_{\Gamma_\ve} |h_\ve|^2\, d\sigma
- \frac{|\Gamma|}{|Y|} \int_\Omega |h_0|^2\, dx\Big]=0.
\end{align*}
Note that the last convergence implies that the two-scale limit $\widetilde{h}_0$ does not depend on $y$. Indeed, by Proposition 2.5 in \cite{AlDa-95}, {one has the estimate}
\begin{align*}
\limsup_{\ve \to 0} \ve \int_{\Gamma_\ve} |h_\ve|^2\, d\sigma\ge
\frac{1}{|Y|} \int_\Omega \int_\Gamma |\widetilde h_0|^2\, d\sigma_y dx
\ge \frac{|\Gamma|}{|Y|} \int_\Omega |h_0|^2 \, dx.
\end{align*}
{Thus, one can see that}
\begin{align*}
\frac{1}{|\Gamma|}\int_\Omega \int_\Gamma |\widetilde h_0|^2\, d\sigma_y dx
= \int_\Omega \left(\frac{1}{|\Gamma|}\int_\Gamma \widetilde h_0\, d\sigma_y\right)^2 \, dx.
\end{align*}
{Moreover, it is clear that}
\begin{align*}
\frac{1}{|\Gamma|}\int_\Omega\int_\Gamma |\widetilde h_0|^2\, d\sigma_y dx
&= \frac{1}{|\Gamma|}\int_\Omega \int_\Gamma |\widetilde h_0 - h_0|^2\, d\sigma_y dx\\
&+ \frac{2}{|\Gamma|}\int_\Omega \int_\Gamma (\widetilde h_0-h_0) h_0 \, d\sigma_y dx\\
&+ \frac{1}{|\Gamma|}\int_\Omega \int_\Gamma |h_0|^2\, d\sigma_y dx = \int_\Omega |h_0|^2\, dx,
\end{align*}
which yields
\begin{align*}
\frac{1}{|\Gamma|}\int_\Omega \int_\Gamma |\widetilde h_0-h_0|^2\, d\sigma_y dx=0 \quad \Rightarrow \quad 
\widetilde h_0 = h_0(t,x).
\end{align*}
Now, {dividing} \eqref{eq:osnovopolagaushaya} by $\delta\not=0$ and passing to the limit as $\delta \to +0$ and $\delta \to -0$, we derive
\begin{align*}
&\frac{|\Gamma|}{|Y|}\int_0^t \int_\Omega\partial_\tau w_0 [\psi]\, dx d\tau
+
\frac{|\Gamma|}{|Y|}\int_0^t\int_\Omega\partial_\tau h_0 \, \varphi\, dx d\tau \nonumber\\
&
+\frac{1}{c_m|Y|}\int_0^t\int_\Omega \int_{Y_i} 
a_i (\partial_{x_1} z_0^i{\bf e}_1+\nabla_y z_{1}^i)\cdot
( \partial_{x_1} \psi^{i}{\bf e}_1 + \nabla_y \Psi^{i}) dy dx  d\tau  \nonumber
\\
&+ \frac{1}{c_m|Y|}\int_0^t \int_\Omega \int_{Y_e} 
a_e (\nabla z_0^e+\nabla_y z_{1}^e)\cdot
( \nabla \psi^{e}+\nabla_y \Psi^{e}) \,dy dx d\tau  \nonumber\\
&+ (\lambda - \frac{1}{c_m}) \frac{|\Gamma|}{|Y|}\int_0^t \int_\Omega w_0[\psi]\, dx d\tau 
-\frac{|\Gamma|}{|Y|c_m}
\int_0^t \int_\Omega h_0 [\psi]\, dx d\tau 
\\
&+(b+\lambda)\frac{|\Gamma|}{|Y|} 
\int_0^t \int_\Omega h_0 \varphi\, dx d\tau  \nonumber
- \theta\frac{|\Gamma|}{|Y|} \int_0^t\int_\Omega w_0 \varphi\, dx d\tau \nonumber\\
&+\frac{|\Gamma|}{3c_m |Y|} \int_0^t\int_\Omega e^{2\lambda \tau} w_0^3 [\psi]\, dx d\tau
-\int_0^t \int_\Sigma \frac{e^{-\lambda \tau}}{c_m} J^e \psi^e\, d\sigma d\tau \nonumber\\
&-  a \frac{|\Gamma|}{|Y|} 
\int_0^t \int_{\Omega} e^{-\lambda \tau}\varphi\, dx d\tau = 0. \nonumber
\end{align*}
Taking $\psi^i=\psi^e=\varphi=0$, we obtain $z_1^e(t, x, y)= N^e(y) \cdot \nabla z_0^e(t,x)$, $z_1^i(t,x,y)= N_1^i(y) \partial_{x_1} z_0^i(t,x)$, where $N_k^e, N_1^i$ solve the cell problems \eqref{eq:cell-prob-ee} and \eqref{eq:cell-prob-ii}, respectively. Note that in the case when $Y_i$ is a cylinder---constant cross-section---, $N_1^i(y)$ is constant. 
Recalling the definition of the effective coefficients $(a_e^\eff)_{kl}$ \eqref{eq:eff-coef-e}, and taking $\Psi^l=0$, 
we obtain 
\begin{align}
&\int_0^t \int_\Omega\partial_\tau w_0 [\psi]\, dx d\tau
+
\int_0^t\int_\Omega\partial_\tau h_0 \, \varphi\, dx d\tau \nonumber\\
&
+\frac{1}{c_m}\int_0^t\int_\Omega  
a_i^\eff \partial_{x_1} z_0^i \,
\partial_{x_1} \psi^{i}dxd\tau 
+ \frac{1}{c_m}\int_0^t \int_\Omega 
a_e^\eff \nabla z_0^e\cdot
\nabla \psi^{e}\,dx d\tau  \nonumber\\
&+ (\lambda - \frac{1}{c_m}) \int_0^t \int_\Omega w_0[\psi]\, dx d\tau 
-\frac{1}{c_m}
\int_0^t \int_\Omega h_0 [\psi]\, dx d\tau \label{eq:weak-effective}\\
&+(b+\lambda)
\int_0^t \int_\Omega h_0 \varphi\, dx d\tau  \nonumber
- \theta \int_0^t\int_\Omega w_0 \varphi\, dx d\tau \nonumber\\
&+\frac{1}{3c_m} \int_0^t\int_\Omega e^{2\lambda \tau} w_0^3 [\psi]\, dx d\tau \nonumber\\
&=\frac{|Y|}{c_m|\Gamma|}\int_0^t \int_\Sigma e^{-\lambda \tau} J^e \psi^e\, d\sigma d\tau + a  
\int_0^t \int_{\Omega} e^{-\lambda \tau}\varphi\, d\sigma d\tau. \nonumber
\end{align}
Performing the change of unknowns $u_0^l= e^{\lambda t} z_0^l$, $v_0=e^{\lambda t} w_0$, $g_0= e^{\lambda t} h_0$, we obtain \eqref{eq:hom-prob}. Note that the initial condition $v_0(0,x)=V^i - V^e$ is obtained using Lemma \ref{lm:compact}. Indeed, for any smooth test function $\psi(t,x)$ such that $\psi(T, x)=0$, we have
\begin{align*}
    \int_\Omega z_\ve^l(0,x) \psi\, dx
    &= \int_0^T \int_{\Omega} (z_\ve^l \partial_t \psi + \partial_t z_\ve^l \psi)\, dx dt \\
    &\mathop{\to}_{\ve \to 0}\quad
    \int_0^T \int_{\Omega} (z_0^l \partial_t \psi + \partial_t z_0^l \psi)\, dx dt= \int_\Omega V^l \psi\, dx.
\end{align*}
The proof of Theorem \ref{th:main-short} is completed.


\section{Well-posedness of the limit problem}
In order to prove the well-posedness of the homogenized problem given by its weak formulation \eqref{eq:weak-effective}, we rewrite it in matrix form as an abstract parabolic equation. We introduce $q_0$ solving the auxiliary problem in $\Omega$:
\begin{align}
\label{eq:q0}
-{\rm div}(a_e^\eff \nabla q_0) &- a_i^\eff \partial_{x_1 x_1}^2 q_0=0, \quad &x&\in \Omega, \nonumber\\
a_e^\eff \nabla q_0\cdot \nu &= \frac{|Y|}{|\Gamma|} J^e, \quad &x& \in \Sigma,\\
q_0&=0, \quad &x&\in S_0\cup S_L.\nonumber
\end{align}
Here, the effective coefficient $a_i^\eff=|Y_i|a_i/|\Gamma|$. 
Multiplication \eqref{eq:q0} by a smooth test function $\psi^e$ such that $\psi^e=0$ on $S_0\cup S_L$ {leads to}
\begin{align}
    \label{eq:int-J}
    \frac{|Y|}{|\Gamma|}\int_\Sigma  J^e \psi^e\, d\sigma= 
    \int_\Omega a_e^\eff \nabla q_0 \cdot \nabla \psi^e\, dx
    + \int_\Omega a_i^\eff \partial_{x_1}q_0 \partial_{x_1} \psi^e\, dx.
\end{align}
Substituting \eqref{eq:int-J} into \eqref{eq:weak-effective}, and introducing
$\widetilde z_0^l=z_0^l - q_0 e^{-\lambda t}$, $l=i, e$, we have the following weak formulation:
\begin{align}
&\int_0^t \int_\Omega\partial_\tau w_0 [\psi]\, dx d\tau
+
\int_0^t\int_\Omega\partial_\tau h_0 \, \varphi\, dx d\tau \nonumber\\
&
+\frac{1}{c_m}\int_0^t\int_\Omega
a_i^\eff \partial_{x_1} \widetilde z_0^i \,
\partial_{x_1} \psi^{i}dxd\tau 
+ \frac{1}{c_m}\int_0^t \int_\Omega 
a_e^\eff \nabla \widetilde z_0^e\cdot
\nabla \psi^{e}\,dx d\tau  \nonumber\\
&+ \left(\lambda - \frac{1}{c_m}\right) \int_0^t \int_\Omega w_0[\psi]\, dx d\tau 
-\frac{1}{c_m}
\int_0^t \int_\Omega h_0 [\psi]\, dx d\tau \label{eq:weak-effective-tilde}\\
&+(b+\lambda)
\int_0^t \int_\Omega h_0 \varphi\, dx d\tau  \nonumber
- \theta \int_0^t\int_\Omega w_0 \varphi\, dx d\tau \nonumber\\
&+\frac{1}{3c_m} \int_0^t\int_\Omega e^{2\lambda \tau} w_0^3 [\psi]\, dx d\tau \nonumber\\
&=a  
\int_0^t \int_{\Omega} e^{-\lambda \tau}\varphi\, d\sigma d\tau
+ \int_0^t \int_{\Omega} e^{-\lambda \tau} a_i^\eff \partial_{x_1 x_1}^2 q_0\, [\psi]\, dx d\tau. \nonumber
\end{align}
We {seek} to rewrite the weak formulation \eqref{eq:weak-effective-tilde} in matrix form as an abstract parabolic equation. To this end, we first introduce the following functional spaces:
\begin{align*}
H_0 &= L^2(\Omega) \times L^2(\Omega), \nonumber\\
H_i& = \{z^i\in L^2(\Omega): \,\,\partial_{x_1} z^i \in L^2(\Omega), \,\, z^i=0 \,\, \mbox{on} \,\, S_0\cup S_L\},\\
H_e& = \{z^e\in L^2(\Omega): \,\,  \nabla z^e \in L^2(\Omega)^3, \,\, z^e=0 \,\, \mbox{on} \,\, S_0\cup S_L\},\nonumber\\
X_0& = \{w=z^i - z^e:\,\, z^i\in H_i, \,\, z^e \in H_e\}.\nonumber
\end{align*}
The norm in $H_i$ is given by
\begin{align*}
\|z\|_{H_i}^2 = \int_\Omega |z|^2 \, dx
+ \int_\Omega |\partial_{x_1} z|^2 \, dx.
\end{align*}
{For the one associated to $H_e$,} we adopt the standard $H^1$-norm.
For each element $w_0\in X_0$, we associate a unique pair $(\widetilde z_0^i, \widetilde z_0^e) \in H_i \times H_e$ solving the following problem
\begin{align}
\label{eq:pair}
-a_i^\eff \partial_{x_1 x_1}^2 \widetilde z_0^i &= {\rm div}(a_e^\eff \nabla \widetilde z_0^e), \quad &x&\in \Omega,\nonumber\\
\widetilde z_0^i - \widetilde z_0^e &= w_0, \quad &x&\in \Omega,\\
a_e^\eff \nabla \widetilde z_0^e \cdot \nu &=0, \quad &x&\in \Sigma,\nonumber\\
\widetilde z_0^i=\widetilde z_0^e&=0, \quad &x&\in S_0\cup S_L.\nonumber
\end{align}
The pair $(\widetilde z_0^i, \widetilde z_0^e)$ can be determined by solving the minimization problem
\begin{align*}
\|w_0\|_{W_0}^2:=\inf\Big\{ 
\int_{\Omega} a_i^\eff |\partial_{x_1} \widetilde z_0^i|^2\, dx
+ \int_{\Omega} a_e^\eff \nabla \widetilde z_0^e \cdot \nabla \widetilde z_0^e\, dx \,\, \big| \,\, \widetilde z_0^i \in W_i, \, \widetilde z_0^e \in W_e\Big\}.
\end{align*}
Note that $W_0$ is a Hilbert space with a scalar product {given by}
\begin{align*}
(w_1, w_2)_{W_0} = 
\int_{\Omega} a_i^\eff \partial_{x_1} z_1^i\, \partial_{x_1} z_2^i\, dx
+ \int_{\Omega} a_e^\eff \nabla z_1^e \cdot \nabla z_2^e\, dx,
\end{align*}
where $(z_1^i, z_1^e)$ and $(z_2^i, z_2^e)$ solve \eqref{eq:pair} for $w_1, w_2$ given.
Now \eqref{eq:weak-effective-tilde} is written in the form
\begin{align*}
\partial_t \begin{pmatrix} w_0 \\[6mm] h_0 \end{pmatrix}
&+ 
\begin{pmatrix}
\displaystyle \frac{1}{c_m}A_\eff w_0 + \frac{1}{c_m}\left(\frac{e^{2\lambda t}}{3}w_0^3 - w_0 - h_0\right) + \lambda w_0 \\[6mm]
 (b+\lambda) h_0 -\theta w_0 
\end{pmatrix}
=
e^{-\lambda t} \begin{pmatrix}
a_i^\eff\partial_{x_1 x_1}^2 q_0\\[6mm]
a 
\end{pmatrix},
\end{align*}
where the operator $A_\eff$ defined on smooth functions $w_0$ by
\begin{align*}
(A_\eff w_0, [\psi])_{L^2(\Omega)}:=\frac{1}{c_m}\int_\Omega
a_i^\eff \partial_{x_1} \widetilde z_0^i \,
\partial_{x_1} \psi^{i}dx 
+ \frac{1}{c_m}\int_\Omega 
a_e^\eff \nabla \widetilde z_0^e\cdot
\nabla \psi^{e}\,dx,
\end{align*}
and $(\widetilde z_0^i, \widetilde z_0^e)$ solve \eqref{eq:pair}.
In operator form {one writes}
\begin{align}
\label{eq:abstract-parabol-eff}
&\partial_t W_0 + \mathbb A_0(t, W_0) = F_0(t), \quad (t,x) \in (0,T)\times \Omega,\\
&W_0(0,x)=W_0^0(x), \quad x\in \Omega. \nonumber
\end{align}
{Therein, we have the following operators}
\begin{align*}
    \mathbb{A}_0(t, W_0) &:= B_0^{(1)}(t, W_0) + B_0^{(2)}(t, W_0),\\
    B_0^{(1)}(t, W_0)&:=\begin{pmatrix}
\displaystyle \frac{1}{c_m}A_\eff w_0 +  (\lambda -\frac{1}{c_m})w_0 - \frac{1}{c_m}h_0\\[4mm]
(b+\lambda) h_0 -\theta w_0 
\end{pmatrix},\\
B_0^{(2)}(t, W_0) &:=
\begin{pmatrix}
\displaystyle
\frac{e^{2\lambda t}}{3c_m} w_0^3 \\[4mm] 0
\end{pmatrix},\\
F_0(t)& := e^{-\lambda t} \begin{pmatrix}
\displaystyle
a_i^\eff\partial_{x_1 x_1}^2 q_0 \\[2mm]
a 
\end{pmatrix}.
\end{align*}
Introducing the spaces
\begin{align*}
    H_0&= L^2(\Omega) \times  L^2(\Omega),\\
    V_1& = X_0\times L^2(\Omega), \quad V_1'=X_0'\times L^2(\Omega),\\
    V_2& = L^4(\Omega) \times L^2(\Omega), \quad V_2'= L^{4/3}(\Omega) \times L^2(\Omega),
\end{align*}
we can prove the existence of a unique solution $W_0 \in L^\infty((0,T); H_0) \cap L^2((0,T); V_1) \cap L^4((0,T); V_2)$ to problem \eqref{eq:abstract-parabol-eff}. It follows, as in Section \ref{sec:existence}, from Theorem 1.4 in \cite{lions1969quelques} and Remark 1.8 in Chapter 2. 

\appendix
\section{Monotonicity method}
{The passage} to the limit in the microscopic problem {requires us to} adapt the method of monotone operators due to G.~Minty \cite{minty1962}. The application of the method to problem \eqref{eq:orig-prob} is given 
in Section \ref{Rabota_Dlya_Nastoyaschih_Patsanov}. The proof is quite technical, and in order to extract the main idea of the method we provide its brief description for a model case when the monotone operator is independent of $\ve$. In \cite{Al-1992}, it is shown how to combine the method of monotone operators and the two-scale convergence for a stationary problem.

Let $A$ be a nonlinear continuous monotone operator in a Hilbert space $H$. The scalar product in $H$ will be denoted by $(u, v)$.  We consider a parabolic problem
\begin{align}
    \label{eq:monotone-prob}
    \partial_t u_\ve + A(u_\ve) =f_\ve,\\
    u_\ve\big|_{t=0}=V_\ve^0. \nonumber
\end{align}
Assume that we know that $u_\ve$ converges weakly to $u_0$, $\partial_t u_\ve$ converges weakly to $\partial_t u_0$, and $f_\ve, V_\ve^0$ converge strongly  in $H$ to $f$ and $V^0$, respectively, as $\ve \to 0$. We {aim} to show that $u_0$ satisfies the limit equation $\partial_t u_0 + A(u_0)=f$. Note that, because of the weak convergence, we cannot pass to the limit in the nonlinear term $A(u_\ve)$ directly.

By monotonicity, for any $w_1, w_2 \in D(A)$, {one has}
\begin{align*}
    (A(w_1) - A(w_2), w_1 - w_2) \ge 0.
\end{align*}
Taking $w_1=u_\ve$, $w_2=u_0+\delta \varphi$, with $\delta \in \mathbb R$ and $\varphi \in C^1([0,T]; D(A))$, and using \eqref{eq:monotone-prob}, we get
\begin{align}
    0&\le
    \int_0^t (A(u_\ve) -A(u_0+\delta \varphi), u_\ve-(u_0+\delta \varphi))d\tau. \nonumber\\
    \label{eq:monot-oper-inequal}
    &=
    \int_0^t(f_\ve, u_\ve-(u_0+\delta \varphi))d\tau
    - \int_0^t(\partial_\tau u_\ve, u_\ve)d\tau + \int_0^t(\partial_\tau u_\ve, (u_0+\delta \varphi))d\tau. \\
    &- \int_0^t(A(u_0+\delta \varphi), u_\ve-(u_0+\delta \varphi))d\tau. \nonumber
\end{align}
Integrating by parts, we get
\begin{align*}
    \int_0^t(\partial_\tau u_\ve, u_\ve)d\tau
    = \frac{1}{2}\int_0^t \frac{d}{d\tau}\|u_\ve\|_{H}^2 d\tau
    = \frac{1}{2}\|u_\ve(t, \cdot)\|_{H}^2 - \frac{1}{2}\|V_\ve^0\|_{H}^2.  
\end{align*}
Then inequality \eqref{eq:monot-oper-inequal} transforms into
\begin{align}
&\frac{1}{2}\|u_\ve(t, \cdot)\|_{H}^2 -
\frac{1}{2}\|u_0(t, \cdot)\|_{H}^2 
-\frac{1}{2}\|V_\ve^0\|_{H}^2
+ \frac{1}{2}\|V^0\|_{H}^2 \nonumber\\
\label{eq:strong-con}
&\le
\int_0^t(f_\ve, u_\ve-(u_0+\delta \varphi))d\tau
- \int_0^t(\partial_\tau u_0, u_0)d\tau\\
&+ \int_0^t(\partial_\tau u_\ve, (u_0+\delta \varphi))d\tau
- \int_0^t(A(u_0+\delta \varphi), u_\ve-(u_0+\delta \varphi))d\tau. \nonumber
\end{align}
Passage to the limit, as $\ve \to 0$, in \eqref{eq:strong-con} {yields}
\begin{align*}
&0\le \frac{1}{2}\limsup_{\ve \to 0} \left(
\|u_\ve(t, \cdot)\|_{H}^2 -
\|u_0(t, \cdot)\|_{H}^2\right)\\
&\le
\delta \int_0^t(-f + \partial_\tau u_0 + A(u_0+\delta \varphi), \varphi)d\tau.
\end{align*}
Since the left-hand side is positive and $\delta$ is arbitrary, that {delivers} the strong convergence of $u_\ve$
\begin{align*}
\limsup_{\ve \to 0} \left(
\|u_\ve(t, \cdot)\|_{H}^2 -
\|u_0(t, \cdot)\|_{H}^2\right)=0.
\end{align*}
Furthermore, {one can show that}
\begin{align}
    \label{eq:aux-monotone}
    \int_0^t (\partial_\tau u_0 + A(u_0+\delta \varphi) - f, \delta \varphi) d\tau \ge 0. 
\end{align}
Dividing \eqref{eq:aux-monotone} first by $\delta>0$ and passing to the limit, as $\delta \to 0$, we obtain
\begin{align*}
    \int_0^t (\partial_\tau u_0 + A(u_0) - f, \varphi) d\tau \ge 0.
\end{align*}
Then, dividing \eqref{eq:aux-monotone} by $\delta<0$ and passing to the limit, as $\delta \to 0$, we have the opposite inequality 
\begin{align*}
    \int_0^t (\partial_\tau u_0 + A(u_0) - f, \varphi) d\tau \le 0.
\end{align*}
Thus,
\begin{align*}
    \int_0^t (\partial_\tau u_0 + A(u_0) - f, \varphi) d\tau = 0.
\end{align*}
The last equality holds for an arbitrary $\varphi \in C^1(0,T; D(A))$, so $\partial_t u_0 + A(u_0)=f$. 

This method is used for problem \eqref{eq:abstract-parabol}, where both the domain and the operator $A$ depend on $\ve$, and the test functions have a more complicated two-scale structure.

\section*{Acknowledgments}
This work is supported by Swedish Foundation for International Cooperation in Research and Higher education (STINT) {with Agencia Nacional de Investigación y Desarrollo (ANID), Chile, through project CS2018-7908 (El Nervio – Modeling Of
Ephaptic Coupling Of Myelinated Neurons)} and Wenner-Gren Foundation.

\bibliography{refs-cells}

\begin{thebibliography}{10}
\expandafter\ifx\csname url\endcsname\relax
  \def\url#1{\texttt{#1}}\fi
\expandafter\ifx\csname urlprefix\endcsname\relax\def\urlprefix{URL }\fi
\expandafter\ifx\csname href\endcsname\relax
  \def\href#1#2{#2} \def\path#1{#1}\fi

\bibitem{standring2021gray}
S.~Standring, Gray's anatomy e-book: the anatomical basis of clinical practice,
  Elsevier Health Sciences, 2021.

\bibitem{hodgkin1952}
A.~L. Hodgkin, A.~F. Huxley, A quantitative description of membrane current and
  its application to conduction and excitation in nerve, The Journal of
  physiology 117~(4) (1952) 500.

\bibitem{fitzhugh1955mathematical}
R.~FitzHugh, Mathematical models of threshold phenomena in the nerve membrane,
  The bulletin of mathematical biophysics 17~(4) (1955) 257--278.

\bibitem{jerez2020derivation}
C.~Jerez-Hanckes, I.~Pettersson, V.~Rybalko, Derivation of cable equation by
  multiscale analysis for a model of myelinated axons, Discrete and Continuous
  Dynamical Systems-Series B 25~(3) (2020) 815--839.

\bibitem{jerez2021multiscale}
C.~Jerez-Hanckes, I.~A. Mart{\'\i}nez, I.~Pettersson, V.~Rybalko, Multiscale
  analysis of myelinated axons, in: Emerging Problems in the Homogenization of
  Partial Differential Equations, Springer, 2021, pp. 17--35.

\bibitem{ramon1978ephaptic}
F.~Ramon, J.~W. Moore, Ephaptic transmission in squid giant axons, American
  Journal of Physiology-Cell Physiology 234~(5) (1978) C162--C169.

\bibitem{bokil2001ephaptic}
H.~Bokil, N.~Laaris, K.~Blinder, M.~Ennis, A.~Keller, Ephaptic interactions in
  the mammalian olfactory system, Journal of neuroscience 21~(20) (2001)
  RC173--RC173.

\bibitem{binczak2001ephaptic}
S.~Binczak, J.~Eilbeck, A.~C. Scott, Ephaptic coupling of myelinated nerve
  fibers, Physica D: Nonlinear Phenomena 148~(1-2) (2001) 159--174.

\bibitem{lin2010modeling}
J.~Lin, J.~P. Keener, Modeling electrical activity of myocardial cells
  incorporating the effects of ephaptic coupling, Proceedings of the National
  Academy of Sciences 107~(49) (2010) 20935--20940.

\bibitem{NEK93}
J.~Neu, W.~Krassowska,
  \href{http://europepmc.org/abstract/MED/8243090}{Homogenization of syncytial
  tissues}, Critical reviews in biomedical engineering 21~(2) (1993) 137---199.
\newline\urlprefix\url{http://europepmc.org/abstract/MED/8243090}

\bibitem{franzone2002degenerate}
P.~C. Franzone, G.~Savar{\'e}, Degenerate evolution systems modeling the
  cardiac electric field at micro-and macroscopic level, in: Evolution
  equations, semigroups and functional analysis, Springer, 2002, pp. 49--78.

\bibitem{pennacchio2005}
M.~Pennacchio, G.~Savar{\'e}, P.~C. Franzone, Multiscale modeling for the
  bioelectric activity of the heart, SIAM Journal on Mathematical Analysis
  37~(4) (2005) 1333--1370.

\bibitem{collin2018mathematical}
A.~Collin, S.~Imperiale, Mathematical analysis and 2-scale convergence of a
  heterogeneous microscopic bidomain model, Mathematical Models and Methods in
  Applied Sciences 28~(05) (2018) 979--1035.

\bibitem{BenMroSaaTal2019}
M.~Bendahmane, F.~Mroue, M.~Saad, R.~Talhouk, Unfolding homogenization method
  applied to physiological and phenomenological bidomain models in
  electrocardiology, Nonlinear Analysis: Real World Applications 50 (2019)
  413--447.

\bibitem{GraKar2019}
E.~Grandelius, K.~H. Karlsen, The cardiac bidomain model and homogenization,
  Networks and Heterogeneous Media 14~(1) (2019) 173--204.

\bibitem{AmaAndTim2021}
M.~Amar, D.~Andreucci, C.~Timofte, Homogenization of a modified bidomain model
  involving imperfect transmission, Communications on Pure and Applied Analysis
  (2021) 0.

\bibitem{Veneroni}
M.~Veneroni, \href{http://dx.doi.org/10.1002/mma.740}{Reaction-diffusion
  systems for the microscopic cellular model of the cardiac electric field},
  Math. Methods Appl. Sci. 29~(14) (2006) 1631--1661.
\newblock \href {https://doi.org/10.1002/mma.740} {\path{doi:10.1002/mma.740}}.
\newline\urlprefix\url{http://dx.doi.org/10.1002/mma.740}

\bibitem{bourgault2009existence}
Y.~Bourgault, Y.~Coudiere, C.~Pierre, Existence and uniqueness of the solution
  for the bidomain model used in cardiac electrophysiology, Nonlinear analysis:
  Real world applications 10~(1) (2009) 458--482.

\bibitem{mandonnet2011role}
E.~Mandonnet, O.~Pantz, The role of electrode direction during axonal bipolar
  electrical stimulation: a bidomain computational model study, Acta
  neurochirurgica 153~(12) (2011) 2351--2355.

\bibitem{nagumo1962active}
J.~Nagumo, S.~Arimoto, S.~Yoshizawa, An active pulse transmission line
  simulating nerve axon, Proceedings of the IRE 50~(10) (1962) 2061--2070.

\bibitem{AlDa-95}
G.~Allaire, G.~Allaire, C.~a.~L. Atomique, A.~Damlamian, U.~Hornung, Two-scale
  convergence on periodic surfaces and applications, Mathematical Modelling of
  Flow through Porous Media (1995).
\newblock \href {https://doi.org/10.1.1.54.8661} {\path{doi:10.1.1.54.8661}}.

\bibitem{lions1969quelques}
J.-L. Lions, Quelques m{\'e}thodes de r{\'e}solution de problemes aux limites
  non lin{\'e}aires (1969).

\bibitem{showalter2013monotone}
R.~E. Showalter, Monotone operators in Banach space and nonlinear partial
  differential equations, Vol.~49, American Mathematical Soc., 2013.

\bibitem{acerbi1992extension}
E.~Acerbi, V.~ChiadoPiat, G.~Dal~Maso, D.~Percivale, An extension theorem from
  connected sets, and homogenization in general periodic domains, Nonlinear
  Analysis: Theory, Methods \& Applications 18~(5) (1992) 481--496.

\bibitem{minty1962}
G.~Minty, Monotone (nonlinear) operators in hilbert space, Duke Math. J. 29
  (1962) 341--346.

\bibitem{Al-1992}
G.~Allaire, \href{http://dx.doi.org/10.1137/0523084}{Homogenization and
  two-scale convergence}, SIAM J. Math. Anal. 23~(6) (1992) 1482--1518.
\newblock \href {https://doi.org/10.1137/0523084} {\path{doi:10.1137/0523084}}.
\newline\urlprefix\url{http://dx.doi.org/10.1137/0523084}

\end{thebibliography}

\end{document}